\documentclass[a4paper, 11pt]{article}

\usepackage[margin=2.5cm]{geometry}

\usepackage{natbib}

\usepackage[latin1]{inputenc}
\usepackage{graphicx}
\usepackage{amsmath,amsthm,bm,mathrsfs}
\usepackage{amsfonts, amssymb}
\usepackage{comment} 
\usepackage{wrapfig}
\usepackage{yfonts, kuvio}  
\usepackage{graphicx}
\usepackage{arydshln}           
\usepackage{psfrag}
\usepackage{nicefrac}
\usepackage{floatflt}
\usepackage{amscd}
\usepackage{ifthen}   

\numberwithin{equation}{section}

\theoremstyle{plain}
\newtheorem{proposition}{Proposition}
\newtheorem{theorem}{Theorem}
\theoremstyle{remark}
\newtheorem{remark}{Remark}
\theoremstyle{definition}
\newtheorem{definition}{Definition}

\newcounter{beispiel}
\newcommand\neuesbeispiel[1]{ \ifthenelse{\equal{
#1}{ohneTitel}}{ \textbf{Example~\thebeispiel:} }{\textbf{Example~\thebeispiel: }(#1). }
} 

\newcommand\beispielende{\  \hspace*{0.97\textwidth} $\triangleleft$ \\[\parskip]} 

\makeatletter

\newenvironment{Bsp}[1][ohneTitel]{%
	\stepcounter{beispiel}
	\small
	\def\@currentlabel{\thebeispiel}
	\ifthenelse{\equal{ohneTitel}{#1}}{\vspace{1ex}\noindent \textbf{Example~\thebeispiel:} }{\vspace{1ex}\noindent\textbf{Example~\thebeispiel: }(#1). }
}{\beispielende}

\makeatother

\newcommand\pd[2]{\frac{\partial #1}{\partial #2}}	
\newcommand\tvec[1]{\ensuremath{\partial_{#1}}}	    
\newcommand\TRaum[2]{\ensuremath{T_{#2} #1}}				
\newcommand\CTRaum[2]{\ensuremath{T_{#2}^* #1}}			
\newcommand\T[1]{T}																	

\newcommand\Rg{\text{rank}\ }					
\newcommand\vecf[1]{\boldsymbol{#1}}								
\newcommand\cvecf[1]{\boldsymbol{#1}}								
\renewcommand\vec[1]{{#1}}													
\newcommand\dMan[1]{\ensuremath{\mathcal{#1}}}
\newcommand\basis[1]{\ensuremath{\mathfrak{#1}}}
\newcommand\mc[1]{\ensuremath{\mathcal{#1}}}
\newcommand\ol[1]{\ensuremath{\overline{#1}}}
\makeatletter

\newcommand{\Rmnum}[1]{\expandafter\@slowromancap\romannumeral #1@}
\makeatother

\newcommand\liekl[2]{\left[\,#1,#2\,\right]}						
\newcommand\prol[2]{\text{pr}^{(#2)} #1}	  				

\newcommand{\pushforward}[2][keinIndex]{
	\ifthenelse{\equal{keinIndex}{#1}}{\ensuremath{#2_*}}{\ensuremath{#2_{*,#1}}}
	}
	
\newcommand{\sprod}[2]{\ensuremath\left<\,#1,\, #2\right>}
\newcommand\pullback[1]{\ensuremath #1^*}

\newcommand{\Span}{ \text{\upshape span}}\,
\newcommand{\cartandistribution}{\ensuremath{\mathcal{C}}}

\def\RR{\mathbb R}																	

\newcommand{\proj}[1]{\ensuremath{\text{proj}_{#1}}}	
\newcommand{\id}{\ensuremath{\text{id}}}	
\newcommand{\indexset}[1]{\ensuremath{\mathbb{#1}}}	
\newcommand{\functionspace}[1]{\ensuremath{\mathcal{#1}}}	

\begin{document}

\title{On symmetries of nonlinear systems in state representation and application of invariant feedback design}
\author{ {\sc Carsten Collon}$^1$, {\sc Joachim Rudolph}$^2$\\[4pt]
 $^1$~Younicos AG,\\ 12489 Berlin, Germany \\[4pt]
 $^2$~Saarland University, Chair of Systems Theory and Control Engineering, \\
 66123 Saarbr\"{u}cken, Germany\\[6pt]
\vspace*{6pt}}
\pagestyle{headings}
\markboth{C.~Collon, J.~Rudolph}{\rm On symmetries of nonlinear systems in state representation and application of invariant feedback design}
\maketitle


\begin{abstract}
{ Symmetries of nonlinear control systems in state representation are considered. To this end, a geometric approach to ordinary differential equations is advocated. Invariant feedback laws for systems with Lie symmetries, i.e.\ feedback laws that preserve the symmetry group of a considered plant, can be constructed based on invariants of the considered  group action. Under minor technical assumptions suitable invariant tracking errors can be determined by following a normalization procedure.  The underlying local geometric meaning of this procedure is discussed  and it is shown how it can also be applied in order to derive a local, reduced-order system representation. Further, the idea of controlled symmetries, i.e.\ imposing desired symmetry properties on a given control system by state feedback, is discussed by outlining an exemplary control design for a predator-prey bioreactor.}
%
{Symmetry,  Lie-Bäcklund equivalence, Lie group, invariant feedback control}

\end{abstract}

\section{Introduction}

Symmetries of dynamical systems have been a subject of long standing interest in the treatment of dynamical systems. 
Roughly speaking, a transformation (or a family of transformations) is a symmetry of a dynamical system if the transformation maps solutions to solutions.  The knowledge of the admitted group of symmetries can help to obtain a structural understanding of the dynamics and the underlying problem that has been modeled. For instance in classical mechanics, the existence of symmetries is closely connected with the existence of conserved quantities and system reduction \citep{AM87}.  
In \cite{GM85} the notion of symmetry is adapted to nonlinear control systems in state representation (i.e.\ under-determined ordinary differential equations with input) and structural consequences of Lie symmetries in the light of a possible decomposition into subsystems of lower dimension plus quadrature are discussed. To this end, an assumed special structure of the Lie algebra (namely Abelian, non-Abelian with center) of the symmetry group is exploited leading to more detailed results than for the general Lie algebra case, which, on the other hand, is closely connected to the notion of {controlled invariant distributions} and {partial symmetries} \citep[cf.][]{NvdS1982, NvdS1985}. The interest in the decomposition of a dynamical system lies in the potential simplification of the control problem (i.e.\ dimensional reduction) as the system equations can be split into one part along the group orbits 
and a reduced order system describing the motion in the orbit space ("from orbit to orbit"), where the latter can be regarded as a local reduced order model that is in some sense complete \citep[see][]{ZZ92}. From a structural point of view symmetries -- as they are intrinsic (i.e.\ coordinate independent) properties of the differential equations -- can serve to classify dynamical systems (see for instance \cite{LR2004} for second order control systems). In both cases the complete knowledge of the symmetry algebra is assumed leading to the problem of the computation of Lie symmetries for a given differential equation. For control systems in state representation this has been considered in \cite{KK94,CKK02}. Further, as shown in \cite{SKZ02}, a symmetry-based approach yields results on non-accessibility and non-observability of nonlinear systems described by implicit differential equations (this should be related to the role of controlled invariant distributions for explicit systems, cf.\ \cite{Isi95}).

In general, symmetries are not invariant under feedback, i.e.\ designing feedback laws using the usual design methods such as feedback linearization, integrator backstepping, or sliding mode control can lead to feedback laws that, when applied, break admitted symmetries of the control problem. Hence, the question arises how feedback laws can be designed that are compatible with a known symmetry group. 

This observation motivates the notion of \emph{invariant feedback}, i.e.\ feedback laws which preserve the error dynamics under the action of the admitted symmetry. To this end, an invariant error approach has been introduced in \cite{RR99} and further examined in \cite{RF03,Rud03}. For a discussion of invariant control design by exact input-output linearization see \cite{MRR04}. The extension of the invariant error approach to invariant asymptotic observers can be found in \cite{AR02, BMR08}. 
 Moreover, the availability of inputs allows the introduction of symmetries into a given control problem (certain restrictions naturally apply). Whereas this is quite naturally part of control design e.g.\ when compensating gravity effects to render a problem symmetric, the injection of problem specific symmetries by feedback has been denoted as \emph{controlled symmetry} in  \cite{SB05}. In the present paper this approach is adopted in order to design a feedback law which is invariant w.r.t.\ to certain changes in the realization of the plant model. 

The amount of literature on  symmetries of dynamical systems, especially those carrying the structure of Lie groups, is by far too extensive to be covered here. Details on the application of Lie groups to differential equations can be found in \cite{BK89,Olv93,Olv95, Ibr94} and the references therein.

The intention of the present paper is twofold. On the one hand it aims to advocate a differential geometric approach to symmetries based on jets and prolongations following the school of Vinogradov  \citep{Vin81, Vin84} and the  approach of \cite{Zha92}. This framework allows an intuitive understanding of symmetries directly motivated from classical geometry, i.e.\ as automorphisms of some geometric object. Further, as this framework has been proved useful w.r.t.\ control-related questions such as the understanding of equivalence and differential flatness \citep[cf.][]{FLMR94, FLMR99}, it also allows the identification of similarities between other geometric results.

On the other hand, the existence of symmetries raises the questions of ``suitable" tracking errors in order to preserve this particular structure. Whereas constructive answers to the question of how to design compatible feedbacks w.r.t.\ to known Lie symmetries have been given in \cite{MRR04},  an example is given in the second half of the paper  to point out the potential use of state-symmetries injected by suitable feedback in order to achieve invariance w.r.t.\ certain transformations.


The article is organized as follows. In Section~\ref{sec:geometry} the geometric approach to differential equations and symmetries is outlined, followed by an intuitive definition of symmetries for differential equations given in Section~\ref{sec:symmetries}. By focusing on symmetry transformations that are elements of a connected Lie group acting locally effectively on a manifold, invariants of these symmetries can be constructed by following a normalization procedure. These invariants play a prominent role in the invariant tracking control approach based on so-called $G$-compatible tracking errors, which is  motivated and outlined in Section~\ref{sec:invariant_control_design}. Further, it is shown how the normalization procedure can be applied  to derive a local reduced order realization of a control system in state representation admitting a state symmetry. Finally, Section~\ref{sec:Beispiel} is devoted to a short example illustrating the use of symmetries induced by suitable state feedback allowing a predator-prey bioreactor to be rendered invariant w.r.t.\ different growth kinetic models. 

\section{A geometrical setting for ordinary differential equations}
\label{sec:geometry}

Since the application of differential geometric objects such as manifolds, distributions etc.\ has become fairly standard in the nonlinear systems literature \citep[see e.g.\!][]{NvdS90, Isi95} common objects such as manifolds and vector fields will not be defined. General introductions to finite-dimensional differential geometry can be found in \cite{War83, Boo03}. The presentation of the geometric framework follows the setting of \cite{Zha92} which has been introduced in a control context before, see  \cite{FLMR94}, \cite{FLMR99}, \cite{Pom95}, \cite{BCD+99}, and \cite{SBS07} for applications, for details on the geometry of jet bundles the reader is referred to \cite{Sau89}.

\subsection{An introductive  example}

\begin{figure}
  \centering \small
  \psfrag{phi}{\hspace{1ex}$\varphi$} \psfrag{phidot}{\hspace{0ex}$\dot{\varphi}$}
  \psfrag{phi2dot}{$\ddot{\varphi}$} \psfrag{S}{$\mc{S}$}
  \psfrag{Phi}{$\Phi$}
  \psfrag{E1}{$\sigma$} \psfrag{E2}{$\tilde{\sigma}$}
 	\includegraphics[width=0.5\linewidth]{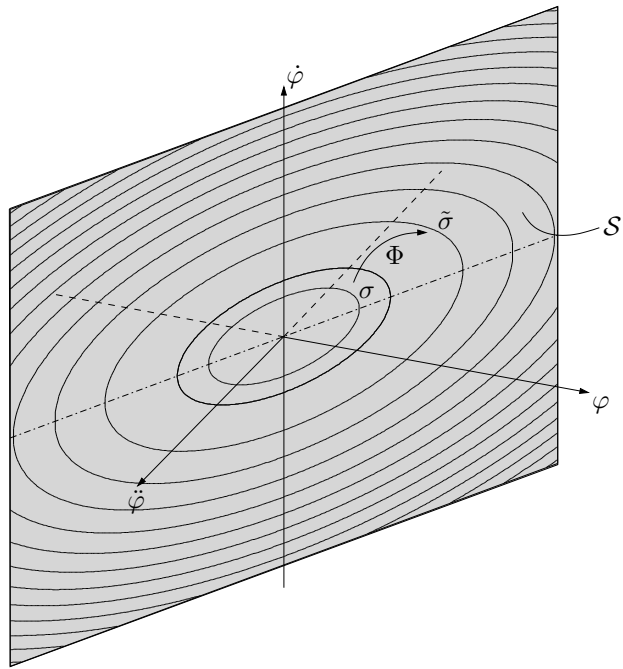} \\
	\caption{Differential equation \eqref{eq:Pendel} as submanifold $\mc{S}$  of the jet space $J^2\pi$ (modulo time reparameterization) and action of $\Phi$ mapping solutions to solutions}
	\label{fig:Pendel_Phasenportrait}
\end{figure}

Consider the equation of a harmonic oscillator with resonance frequency $\omega$ 
\begin{align}
		\ddot{\varphi}(t) &= - \omega^2 \varphi(t),  \qquad t\in I\subset \RR, \, \varphi(t)\in \RR,
		\label{eq:Pendel}
\end{align}   
where the trivial bundle $\left(M, \pi=\proj{1}, I\right)$, $M=I\times \RR$, with global coordinates $(t, \varphi)$ is used to describe the evolution of its configuration w.r.t.\ time. %
Given a smooth function $\sigma: I \rightarrow M$, $t \mapsto (t, \sigma(t))$ (i.e.\ a graph),  $\sigma$ is called a solution of the differential equation if $\sigma$ together with its second time derivative fulfills \eqref{eq:Pendel}. Turning to the geometric description a solution can be interpreted as a graph given by  $\sigma$ and its time derivatives in the extended space $J^2 \pi$, the second jet space,  $j^2\sigma: I \rightarrow J^2\pi$, $t\mapsto \left(t, \sigma(t), \dot{\sigma}(t), \ddot{\sigma}(t)\right)$, where $j^2 \sigma$ is called the second prolongation of $\sigma$. In the jet space with local coordinates $(t, \varphi, \dot{\varphi}, \ddot{\varphi})$ the differential equation~\eqref{eq:Pendel} defines a regular submanifold 
\begin{align} \label{eq:pendel_untermannigfaltigkeit}
			J^2 \pi \supset \mc{S} = \left\{ \left(t, \varphi, \dot{\varphi}, \ddot{\varphi}\right)\in J^2\pi:\, \ddot{\varphi} + \omega^2 \varphi = 0 \right \}
\end{align}
 introducing a geometrical object for the differential equation. As the differential equation is time-invariant, it is sufficient to look at the phase space for visualization purposes obtaining the well-known phase portrait shown in Figure~\ref{fig:Pendel_Phasenportrait}. A function  $\sigma$ is a solution iff the graph of its prolongation \linebreak $t\mapsto (t, \sigma(t), \dot{\sigma}(t), \ddot{\sigma}(t))$ is contained in the regular submanifold:  $j^2 \sigma(t) \subset \mc{S}$, $\forall t\in I$. 

By means of the geometric approach, i.e.\ regarding solutions as graphs of smooth functions on a manifold, the understanding of symmetry of a differential equation can be related to the classical geometric idea of symmetry, i.e.\ an automorphism of a given geometric object is understood as symmetry hereof (e.g.\ the reflection about the center of a circle). Returning to the example, consider the mapping $\Phi: \left(t, \varphi, \dot{\varphi}, \ddot{\varphi}\right) \mapsto \left(t, \lambda \varphi, \lambda  \dot{\varphi}, \lambda \ddot{\varphi}\right)$, $\lambda \in \RR\backslash\{0\}$. Applying $\Phi$ to the definition of $\mc{S}$ given in \eqref{eq:pendel_untermannigfaltigkeit} yields
\begin{align*}
		\Phi \left(S\right) = \left\{ \left(t, \varphi, \dot{\varphi}, \ddot{\varphi}\right)\in J^2\pi,\,: \lambda\left( \ddot{\varphi} + \omega^2 \varphi \right) = 0 \right \} = \mc{S},
\end{align*}
i.e.\ the submanifold is invariant w.r.t.\ $\Phi$. 
Each solution $\sigma$  is mapped to another solution $\tilde{\sigma}$ as $\mc{S}$ is invariant, and due to this characteristic property $\Phi$ is denoted a symmetry of the differential equation.   Of course this can can be observed in coordinates as well: Let $\tilde{\varphi}$, $\dot{\tilde{\varphi}}$, $\ddot{\tilde{\varphi}}$ denote the transformed variables, then the transformed differential equation reads
\begin{align*}
		\lambda \ddot{\varphi}(t) &= \ddot{\tilde{\varphi}} =  - \omega^2 \lambda \varphi = -\omega^2 \tilde{\varphi},
\end{align*}
showing that the differential equation is invariant w.r.t.\ $\Phi$. 

\subsection{Notions from differential geometry}

As sketched in the introductory example the geometric treatment of ordinary differential equations in this paper aims to allow a definition of symmetry that is strongly motivated by the usual geometric understanding of symmetries. Hence a possible identification of a system of ordinary differential equations (in state space representation) with submanifolds of jet manifolds is outlined based on the identification of the differential equation with its totality of solutions which are locally given  as graphs of smooth functions. Before this geometric picture can be established in a more formal way some basic notions from differential geometry are recalled. 

A smooth fibered manifold is a triple $\left(\mc{E}, \pi, \mc{B}\right)$ consisting of the total manifold $\mc{E}$,  the base manifold $\mc{B}$, and a smooth surjective submersion\footnote{Sometimes $\pi$ is referred to as projection which it is only in adapted coordinates.} $\pi: \mc{E} \rightarrow  \mc{B}$. The following considerations are restrained to the finite dimensional case, more specifically $\dim \mc{E} = 1 + q$,  and adapted coordinates on $\mc{E}$ and $\mc{B}$ are given by $(t, z^1, z^2, \hdots, z^q)=\left(t, \vec{z}\right)$ and $(t)$, respectively. 

For each $t\in \mc{B}$ the total manifold is locally diffeomorphic to the product space $\mc{E}\stackrel{\mathrm{loc}}{\simeq} \mc{B} \times \mc{F}_t$,  where $\mc{F}_t$ denotes the fiber $\mc{F}_t = \pi^{-1}(t)$ over $t$. Assuming that all fibers $\mc{F}_t$ are diffeomorphic to a typical fiber $\mc{F}$ one arrives at a fiber bundle  $\left(\mc{E}, \pi, \mc{B}\right)$. In order to simplify notation the bundle is usualle only denoted by $\pi$ instead of $\left(\mc{E}, \pi, \mc{B}\right)$. A section $\sigma$ in $\pi$ is a map $\sigma: \mc{B} \rightarrow \mc{E}$ such that $\pi \circ \sigma = \id$ holds wherever $\sigma$ is defined. The set of all sections in $\pi$ is denoted by $\Gamma(\pi)$,  the set of local sections defined around some $\vec{p}\in \mc{E}$ is denoted by $\Gamma_{\vec{p}}(\pi)$.

In order to be able to discuss differential equations geometrically the time derivatives of the dependent coordinates $z^i$ have to be included within the geometric picture leading to the manifold of $k$-jets. Two local sections $\sigma, \psi \in \Gamma_{\vec{p}}(\pi)$ with $\sigma(\vec{p})=\psi(\vec{p})$ and further coinciding in their derivatives up to the $k$-th order are called $k$-equivalent in $\vec{p}$. The equivalence class formed by all sections being $k$-equivalent in $\vec{p}$ to which $\sigma$ belongs is called the $k$-jet of $\sigma$ in $\vec{p}$ and is denoted by  $j^k_{\vec{p}} \sigma$. The reader can think of a $k$-jet locally as the equivalence class of Taylor series expansions of smooth functions that coincide up to the $k$-th order. The set of all $k$-jets for all $\sigma\in \Gamma(\pi)$ forms the $k$-th jet manifold $J^k \pi$. Adapted coordinates $(t,\vec{z})$ induce adapted coordinates on $J^k\pi$ given by $(t, \vec{z}^i, z^i_j)$, $i=1,\hdots,q$, $j=1,\hdots, k$. Together with $\mc{B}$ one defines the bundle $\left(J^k\pi, \pi_k, \mc{B}\right)$. A section $\sigma\in\Gamma_{\vec{p}}(\pi)$ induces a mapping on $J^k \pi$ which is given by its prolongation $j^k \sigma: t\mapsto (t, \sigma(t), \dot{\sigma}(t), \hdots, \sigma^{(k)}(t))$. Given a local  section $\psi\in \Gamma_{\vec{p}}(\pi_k)$, $\psi: \mc{B} \rightarrow J^k\pi$, $t\mapsto \left(t, \psi(t), \psi_1(t), \hdots, \psi_k(t)\right)$ does not necessarily have to be the prolongation of any section in $\pi$. Since solutions of differential equations are prolongations of sections in $\pi$ the contact distribution
\begin{align} \label{eq:kontaktdistribution}
			\cartandistribution^k: \cvecf{\omega}_j^i = dz^i_j - z^i_{j+1} dt = 0, \qquad i= 1, \hdots, q, \, j=0,1,\hdots, k-1,
\end{align}
defined via $(qk)$ so-called contact forms, is used to distinguish arbitrary sections in $J^k \pi$ from prolongations of sections in $\pi$. A local section $\psi$ in $J^k \pi$ is the prolongation of some section $\phi$ in $\pi$ iff the contact conditions~\eqref{eq:kontaktdistribution} hold for $\psi$. 
A smooth mapping $\Phi:J^k \pi  \rightarrow  J^k \pi$ that preserves the contact distribution, i.e.\ $\pushforward{\Phi}\!\left(\cartandistribution^k \right)\subset \cartandistribution^k$, is called contact transformation, i.e.\ $\Phi$ maps prolongations of smooth sections (graphs) to prolongations of (other) smooth sections, where $\pushforward{\Phi}$ denotes the push-forward along $\Phi$.  For more details on the geometry of jet spaces the reader is referred  to \cite{Sau89}.

\subsubsection{Differential manifold and Lie-Bäcklund mappings} 

In the following a differential manifold is distinguished from the usual notion of smooth manifolds by defining a differential manifold  to be a pair $\dMan{M} = (M, \cartandistribution_M)$ consisting of a smooth manifold $M$ and an involutive Cartan distribution 
 $\cartandistribution_M$ \citep[for details cf.][]{Zha92}. The Cartan distribution is finite dimensional and its dimension $m$ is called the Cartan dimension of the differential manifold. The underlying manifold $M$ can be in general infinite dimensional with coordinates indexed w.r.t.\ a countable index set $\indexset{A}$ with cardinality $\#\indexset{A}$. If not stated otherwise it is assumed that each index set is finite dimensional. Note that every finite-dimensional smooth  manifold can be made into a differential manifold in this sense by assigning the tangent bundle as the Cartan distribution, i.e.\ $\dMan{M}=(M, \TRaum{M}{})$. 
A (standard) chart on $\dMan{M}$ is given by a quadruple $(U, \varphi, \RR^{m+\indexset{A}}, \basis{d}_M)$ consisting of a chart $(U, \varphi)$ for $M$ and a basis $\basis{d}_M$ of the Cartan distribution spanned by $m$ base fields.  Since the Cartan distribution is involutive for each point $\vec{p}\in M$ there exists a unique integral manifold of dimension $m$ passing through $\vec{p}$. 
The standard charts on $\dMan{M}$ introduce a natural separation between the $m$ independent coordinates (usually denoted by $\vec{x}^i$, $i=1,\hdots, m$) and the dependent coordinates (usually denoted $u^{\alpha}$, $\alpha\in\indexset{A}$) together forming the adapted coordinates $\left(x^i, u^{\alpha}\right)$ on $M$. Integral manifolds of $\cartandistribution_M$ will be of particular interest as they will play the role of solutions of a system of differential equations. Every integral submanifold $S\subset \dMan{M}$ can be locally written in parametric form $\dMan{M} \supset S : u^{\alpha} = s^{\alpha}(\vec{x})$,  $ \alpha \in\indexset{A}$, i.e.\ as a graph of smooth functions $s^{\alpha}$.  In the following all considerations are restricted to the case $m=1$ (i.e.\ ordinary differential equations) and $\#\indexset{A} = q$.  The independent coordinate will be denoted $t$ (time) and in order to distinguish from the general case the adapted coordinates are denoted  by $\left(t, z^1, z^2, \hdots, z^q\right)=\left(t,\vec{z}\right)$.
The set of smooth functions on $M$ is ${C}^{\infty}(M)$  
and the set of vector fields on $M$ with the same domain $\mathcal{O}$ will be denoted by $\functionspace{T}(\mathcal{O})$.

The motivation for the distinction between smooth and differential manifolds by introducing the Cartan distribution can be understood from the leading example: Let us regard the plane $\mc{S}$ in Figure~\ref{fig:Pendel_Phasenportrait} as a geometric object that is to be identified with the differential equation as a subset of $J^1 \pi$ with adapted coordinates $(t, \varphi, \dot{\varphi})$ rather than in $J^2 \pi$, i.e.\ elements of $\mc{S}$ that are prolongations of graphs that fulfill the differential equation. In general, $\mc{S}$ also contains graphs of prolongations that do not comply with the differential equation. The distinguishing element between solutions and non-solutions is their tangent space in each point. Whereas for solutions the tangent space is uniquely determined by the differential equation, it is not restricted for general graphs. If in the leading example the Cartan distribution for $M=J^1\pi$ is defined as $\cartandistribution(M)= \Span\left\{ \tvec{t} - \omega^2 \varphi \tvec{\dot{\varphi}} \right \}$ a differential manifold $\left(\mc{S}, \cartandistribution(M)\right)$ is obtained in the above sense consisting only of prolonged graphs fulfilling the differential equation. 

Now, consider a smooth mapping $\Phi:M \rightarrow N$ between two differential manifolds  $\dMan{M}=(M, \cartandistribution_M)$ and $\dMan{N}=(N, \cartandistribution_N)$. The mapping is called a Lie-Bäcklund mapping if it preserves the Cartan distribution, i.e.\ if $\pushforward{\Phi}(\cartandistribution_M) \subset \cartandistribution_N$ holds. 
Regarding the coordinate form of a Lie-Bäcklund mapping the following result is recalled.

\begin{proposition}[\citealp{Zha92}] Let $\dMan{M}=(M, \cartandistribution_M)$ and $\dMan{N}=(N, \cartandistribution_N)$  be differential manifolds of Cartan dimension one. Let $(U, \varphi, \RR\times \RR^{\indexset{A}}, \basis{d}_t)$ and $(N, \psi,  \RR\times \RR^{\indexset{B}}, \basis{d}_{\tilde{t}})$ on $\dMan{M}$ and $\dMan{N}$ denote their standard charts with coordinates $(t, \vec{u})$,  $(\tilde{t},\vec{v})$  and let their Cartan distributions be spanned by the two vector fields
\begin{align*}
	\vecf{\partial}_{t} = \tvec{t} + \sum_{\alpha\in\indexset{A}} A^{\alpha}({t},\vec{u}) \tvec{u^{\alpha}}, \qquad
	\vecf{\partial}_{\tilde{t}} = \tvec{\tilde{t}} + \sum_{\beta\in\indexset{B}} B^{\beta}(\tilde{t},\vec{v}) \tvec{v^{\beta}}.
\end{align*}
A smooth mapping $f:U \rightarrow V$ given in coordinates by 
\begin{align}
  \label{eq:LieBaecklundAbbildung_Form}
	\tilde{t} = T({t}, \vec{u}),\quad  v^{\beta} = V^{\beta}({t},\vec{u}), \quad \beta\in\indexset{B}
\end{align}
is a Lie-Bäcklund mapping iff the coordinate functions fulfill the following (defining) equations
 \begin{align}
 	\vecf{\partial}_{t} (V^{\beta}) -  B^{\beta}(T, \vec{V}) \vecf{\partial}_{t}(T) = 0, \quad \beta\in\indexset{B}.
 	\label{eq:LieBaecklund_definierendeGln}
 \end{align}
\end{proposition}

In the following, tangency of vector fields w.r.t.\ submanifolds will be of importance. 

\begin{proposition}[cf.\ \citealp{Zha92}]\label{beh:tangentielles_vektorfeld} A vector field $\vecf{v}\in\functionspace{T}(M)$ is tangent to a submanifold  $S\subset M$, i.e.\ $\left.\vecf{v}\right|_S \in \functionspace{T}(S)$, iff $\vecf{v}:\functionspace{T}(S)\rightarrow \functionspace{T}(S)$, that is iff $			\left.\vecf{v}(s)\right|_S = 0\, \text{for } \left.s\right|_S=0.$
\label{beh:Teilmannigfaltigkeit_tangentialesVektorfeld}
\end{proposition}

A mapping $f$ is said to be a Lie-B\"{a}cklund isomorphism from $\dMan{M}$ to $\dMan{N}$ at the pair of points $(\vec{p}, \vec{q})$ with $\vec{p}\in M$ and $\vec{q}\in N$, if $f$ is a Lie-B\"{a}cklund mapping with $\pushforward{f}(\cartandistribution_M) = \cartandistribution_N$, and if $f$   has a smooth inverse $g$ from a neighborhood of $\vec{q}$ in $N$ to a neighborhood of $\vec{p}$ in $M$  with $\pushforward{g}(\cartandistribution_N) = \cartandistribution_M$. If this is the case, the differential manifolds $\dMan{M}$ and $\dMan{N}$ are said to be (differentially) equivalent. 
If $\Phi:\dMan{M}\rightarrow\dMan{N}$ is a Lie-B\"{a}cklund isomorphism between two differential manifolds $\dMan{M}$ and $\dMan{N}$ with smooth inverse $\Psi:\dMan{N}\rightarrow \dMan{M}$, and $S: \dMan{N} \rightarrow \dMan{N}$ is an automorphism on $\dMan{N}$, i.e. $\pushforward{S}(\cartandistribution(\dMan{N}))\subset  \cartandistribution(\dMan{N})$, then the mapping $\Psi \circ S \circ \Psi^{-1}:\dMan{M} \rightarrow \dMan{M}$ defines an automorphism on $\dMan{M}$.

\subsection{Geometric interpretation of control systems in classical state representation}

Consider the smooth manifolds $\mc{B}$, $\mc{E}$, and $\mc{U}$ referred to as base, state, and input manifold with adapted coordinates $(t)$, $\left(t, x^1, x^2, \hdots, x^n\right)=\left(t,\vec{x}\right)$, and $\left(t, u^1,\hdots, u^m\right)=(t,\vec{u})$, respectively. From these manifolds the bundles $\left(\mc{E}, \pi, \mc{B}\right)$ and $\left(\mc{U}, \pi_u, \mc{B} \right)$ are defined. Using the bundle product over the same base manifold $\times_{\mc{B}}$  yields the manifold $M=\mc{E} \times_{\mc{B}} \mc{U}$ being the total manifold of the bundle $\left(M, \pi_M = \pi \times_{\mc{B}} \pi_{u}, \mc{B} \right)$ with adapted coordinates $\left(t, \vec{x}, \vec{u}\right)$. The mapping $\rho:=\pi^*\left(\pi_{u}\right): M \rightarrow \mc{E}$ can be used to pull back the tangent bundle $\left(\TRaum{\mc{E}}{}, \tau_{\mc{E}}, \mc{E}\right)$ to $M$ yielding the pull-back bundle $\left(\pullback{\rho}(\TRaum{\mc{E}}{}), \pullback{\rho}(\tau_{\mc{E}}), M\right)$ with  $\rho^*\left(\TRaum{\mc{E}}{}\right)\simeq \left(t,\vec{x}, \vec{u}, \dot{t}, \dot{\vec{x}}\right)$. The following commutative diagram summarizes the bundle construction.
\[
 \Diagram 
  &&&	 \rho^*\left(\TRaum{\mathcal{E}}{}\right) &&\lTo^{\rho^*} &\TRaum{\mathcal{E}}{}  \\
  &&&	 \dTo <{\rho^*\left(\tau_{\mathcal{E}}\right)} \rt{-1pt}  \uTo >{f} \rt{3pt} 
       &&& \dTo >{\tau_{\mathcal{E}}}  \\ 
  {\mc{U}} & \lTo^{\pi_{u}^*\left(\pi\right)} && \mc{E} \times_{\mathcal{B}} \mathcal{U} && \rTo^{\rho:=\pi^*\left(\pi_{\vec{\alpha}}\right)} \up{3pt} 
\lTo_{\sigma} \up{-2pt}   &\mathcal{E} \\
   \dTo<{\pi_{u}} &&& \dTo >{\pi \times_{\mathcal{B}} \pi_{u} }  &&& \dTo >{\pi}\\
   \mathcal{B} & \rTo^{\mathrm{id}_{\mathcal{B}}} && \mathcal{B} && \rTo^{\mathrm{id}_{\mathcal{B}}} & \mathcal{B}\\
  \endDiagram
\]

Now, consider a system of ordinary differential equations in classical state representation on $M$, i.e.\
\begin{align} 
	\dot{\vec{x}}^i(t) &= \frac{d}{dt} \vec{x}^i(t) = {f}^i\left(t, \vec{x}(t),\vec{u}(t)\right), \qquad i=1,2,\hdots, n. 
	\label{eq:zustandssystem2} 	
\end{align}
The smooth functions $f^i\in {C}^{\infty}(M)$ define a vector field\footnote{Throughout the paper summation over repeated indices is assumed.}
\begin{align} \label{eq:vektorfeld}
				\Gamma\left(\rho^*\left(\TRaum{\mc{E}}{} \right)\right)\owns \vecf{v}_f &= \tvec{t} + f^i\left(t, \vec{x}, \vec{u}\right) \tvec{x^i}
\end{align}
which is a section in the pull-back bundle $\rho^*\left(\TRaum{\mc{E}}{} \right)$. A local section $\sigma: \mc{B} \rightarrow M$, $t\mapsto (t, \vec{x}(t), \vec{u}(t)) \in \pi_M^{-1}(t)$ whose first prolongation fulfills the system equations~\eqref{eq:zustandssystem2} is called a solution of the system. Since $f=(f^i)$ depends on the input $\vec{u}$ the system is under-determined. 
One arrives at the following two possible geometric objects that can be understood as representatives of the differential equation~\eqref{eq:zustandssystem2}: 
\begin{itemize}
\item \emph{Differential manifold perspective:} The vector field~\eqref{eq:vektorfeld} spans the Cartan distribution $\cartandistribution_f=\Span\{\vecf{v}_f \}$ of the differential manifold $\dMan{M}_f=\left(M, \cartandistribution_f \subset \pullback{\rho}\left(\TRaum{\mc{E}}{}\right)  \right)$, 
where the Cartan distribution  distinguishes prolongations of solutions from prolongations of other sections.  The geometric object identified with the differential equation is given by the differential manifold $\dMan{M}_f$ above.
\item \emph{Jet manifold perspective:}
The first jet manifold $J^1 \pi_M$ has adapted coordinates $\left(t, \vec{x}, \dot{\vec{x}}, \dot{\vec{u}}\right)$, and herein the differential equation \eqref{eq:zustandssystem2} defines a regular submanifold $\mathcal{S}\subset J^1 \pi_M$ given by
\begin{align} \label{eq:submanifold_defining_equations}
	\mathcal{S}=\left\{ \vec{p} \in J^1 \pi_M \, | \, F(\vec{p})=0 \right\}\quad \text{with } F(\vec{p})= \dot{\vec{x}} - f\left(t, \vec{x}, \vec{u}\right).
\end{align} 
 Each smooth section $\sigma:t \mapsto \left(\vec{x}(t), \vec{u}(t)\right)$  has a prolongation $j^1 \sigma$ to $J^{1}\pi_M$ of the form\linebreak $t\mapsto \left(t, \vec{x}(t), \vec{u}(t), \dot{\vec{x}}(t),\dot{\vec{u}}(t)\right)$, where $\sigma$ is a solution of \eqref{eq:zustandssystem2} iff $j^1 \sigma \subset \mathcal{S}$ holds. As  geometric object  
 the embedded closed submanifold $\mathcal{S}\subset J^1 \pi_M$ is identified with the system~\eqref{eq:zustandssystem2}. 
\end{itemize}

\begin{remark} \label{rem:contact_cartan}
	As mentioned before, as long as the tangent bundle is not restricted to the contact distribution, the submanifold $\mathcal{S}$ also contains sections that are not prolongations of a section in $\pi_M$. Consequently, it is natural to equip $\mathcal{S}$ with a Cartan distribution being exactly determined by the contact forms \citep{BCD+99}, thus, in the following $\mathcal{S}$ stands for the pair $(\mathcal{S}, \TRaum{\mathcal{S}}{} \cap \functionspace{C}^1)$. However, since in this case the Cartan distribution is not involutive, the pair does not define a differential manifold in the sense of the definition above, even though the construction serves a similar purpose.
\end{remark}

\subsubsection{Equivalence under Lie-B\"{a}cklund mappings}

By identifying differential equations with differential manifolds one obtains a notion of equivalence of systems under Lie-B\"{a}cklund mappings. The following notion of equivalence has been already discussed in a state space context in \cite{FLMR94, FLMR99}. 
Two systems $\dMan{M}=(M, \Span\{\vecf{v}_f\})$ and $\dMan{N}=(N, \Span\{\vecf{v}_g\})$ are (orbitally) equivalent at $(\vec{p}, \vec{q})\in M\times N$, iff there exists a Lie-B\"{a}cklund isomorphism $\Phi$ between open neighborhoods of $\vec{p}$ and $\vec{q}$ such that $\vec{q}=\Phi(\vec{p})$. Since both differential manifolds consist of integral curves of their Cartan distribution the trajectories of both systems are $\Phi$-related in $(\vec{p},\vec{q})$.  The systems are  (orbitally) equivalent if there exists such a mapping from an open dense subset $U\subset M$ to $N$. The systems are (orbitally) equivalent if such a mapping exists for any pair $(\vec{p}, \Phi(\vec{p}))$ on a open dense subset of $M$. 

Clearly, $\Phi$ maps $\vecf{v}_f$ to an element $\phi\cdot \vecf{v}_g$, $\phi\in {C}^{\infty}(N)$, of $\Span\{\vecf{v}_g\}$. If the parametrization of the integral curves is preserved, i.e.\ $\phi\equiv 1$, the two systems are called differentially equivalent. The latter case is related to the usual notion of $\phi$-related flows of vector fields in the finite-dimensional case \citep{War83}.

As shown in \cite{FLMR99} the input dimension $m$ is preserved under Lie-B\"{a}cklund isomorphisms but the state dimension $n$ can be lost. Hence such equivalences can be useful in order to simplify control related problems, i.e.\ feedback design or the search for symmetries by reduction of the state dimension. 

\begin{remark}
  It should be pointed out that it can be useful to look at a system representation of lower state dimension that is Lie-Bäcklund equivalent to the original system representation in order to simplify the computation of admitted symmetries. For instance, consider the equations
\begin{align*}
	\ddot{y}^1 &= -u^1 \sin \theta + \epsilon u^2\cos\theta, & 	\ddot{y}^2 &= \phantom{-}u^1  \cos \theta + \epsilon u^2\sin\theta  +g, & 	\ddot{\theta} &= u^2,
\end{align*}
describing the planar motion of a rigid body modelling a  PVTOL vehicle (cf.\ \citealp{HSM92}). Here $\left(y^1, y^2\right)$ is the position of the center of mass, $\theta$ describes the orientation of the body, $g$ denotes the gravitational acceleration, the inputs $u^1$ and $u^2$ are the linear and the rotational acceleration, and $\epsilon\ll 1$ is a parameter describing the coupling of the equations w.r.t.\ $u^2$. This system is Lie-Bäcklund equivalent to 
\begin{align*}
	\ddot{z}^1 &= -v^2 \sin v^1, && \ddot{z}^2 = v^2 \cos v^1 + g,	
\end{align*}
with $z^1=y^1-\epsilon \sin\theta$, $v^1=\theta$, $z^2=y^2+\epsilon\cos\theta$,  $v^2 = u^1 - \epsilon \dot{\theta}^2$. Any symmetry of the latter system is reflected by a symmetry of the original equations. Hence in terms of the computation of symmetries it might be easier to consider the representation in $(z,v)$-coordinates. 
\end{remark}

\begin{remark} \label{rem:flatness}
	In \cite{FLMR99} the notion of differential flatness of nonlinear control systems is presented within a differential geometric framework based on Lie-Bäcklund equivalence. To this end, consider an infinite-dimensional manifold $\RR_m^{\infty}$ with coordinates $\left(y^i_{\nu}\right)$, $i=1,\hdots, m$, $\nu\geq 0$, together with the Cartan distribution $\cartandistribution=\Span\left\{ \sum_{\nu\geq 0} \, y^i_{\nu+1} \tvec{y^i_{\nu}}\right\}$, where $y^i_{\nu}$ denotes the $\nu$-th time derivative coordinate of the $i$-th component of $\vec{y}$, i.e.\ a chain of $m$ independent integrators of arbitrary length. A control system $\dMan{M}_f=\left( M, \Span\{ \vecf{v}_f\}\right)$ is differentially flat iff it is Lie-Bäcklund equivalent to a trivial system. 
\end{remark}

\section{Classical symmetries and symmetry groups} 
\label{sec:symmetries}

As motivated in the leading example the differential geometric framework the definition of symmetry for differential equations can be understood from the classical geometrical meaning.

\begin{definition} A symmetry of a differential equation $\mathcal{S}\subset J^1 \pi_M$ is an automorphism on $\mathcal{S}$.
\end{definition}
As pointed out in Remark~\ref{rem:contact_cartan} it is assumed that $\mathcal{S}$ consists of integral curves of the contact distribution. Therefore, any symmetry transformation is necessarily a contact transformation on $J^1\pi_M$ which is equivalent to a Lie-Bäcklund transformation on $\left(\mathcal{S}, \TRaum{\mathcal{S}}{} \cap \functionspace{C}^1 \right)$. 

\begin{definition}	A contact transformation  $\Phi: J^1\pi_M \rightarrow J^1 \pi_M$ for which $\Phi(\mathcal{S})\subset\mathcal{S}$  holds is called a \emph{classical symmetry} of the differential equation $\mathcal{S}$. A smooth vector field $\vecf{v}_{\Phi}\in \Gamma(\TRaum{J^1\pi}{})$ is called an \emph{infinitesimal symmetry} of $\mathcal{S}$ if its flow is a classical symmetry.
\end{definition}

From this definition it is not clear how a contact transformation on $J^1 \pi_M$ is related to point transformations on $M$, which have been initially considered. In the finite-dimensional case this question is answered by a result given by \cite{Bkl75}. For details see also \cite{AI79, BCD+99}.

\clearpage
\begin{theorem}[\citealp{Bkl75}] \label{thm:baecklund}
	Any contact transformation on the jet manifold $J^k \pi_M$, $k\geq 1$, is
\begin{enumerate}
	\item $\dim \pi=q=1$: the $(k-1)$-th prolongation of some contact transformation on $J^1 \pi_M$,
	\item $\dim \pi=q>1$: the $k$-th prolongation of some point transformation (diffeomorphism) of the independent and independent variables on $J^0 \pi_M = \pi_M$.
\end{enumerate}	
\end{theorem}
\noindent Therefore, as $n+m$ will always exceed one, all symmetries derive as prolongations of diffeomorphisms on $\pi_M$.

\begin{remark} 
	Following the nomenclature of \cite{BCD+99} the prefix \emph{classical} is used for symmetries that preserve the order of the differential equation (i.e.\ the considered transformations close off w.r.t.\ the order of the highest derivatives).  The mentioned result of Bäcklund is exactly due to this restriction. Since this type of symmetries is defined independently from the differential equation, they are referred to as \emph{external} -- they are also defined ``outside" of $\mathcal{S}$. Symmetry transformations that are only defined on $\mathcal{S}$ are consequently denoted \emph{internal}, and any external symmetry restricted to $\mathcal{S}$ generates an internal symmetry. However, if one allows the considered transformations only to be declared on the domain of the differential equation, i.e.\ one can use the differential equation in order to re-express higher-order time derivatives occurring due to transformation that do not preserve the order, one arrives at so-called \emph{dynamical symmetries}, which are by definition internal \citep{OAK93}. More generally, if an infinite-dimensional framework is used (i.e.\ for partial differential equations), Bäcklund's result no longer holds allowing the introduction of \emph{generalized symmetries} \citep{Vin84, BCD+99, Olv93}, i.e.\ Lie-Bäcklund transformations which do not preserve the order of the differential equation. For ordinary differential equations generalized symmetry can be related to a dynamical symmetries, see \citep{OAK93} for details.
\end{remark}

\begin{Bsp}
	\label{bsp:LieBaecklundForminvarianz}
	Consider a system of ordinary differential equations in state space representation 
	\begin{align*}
			\dot{\vec{z}}^i &= f^i(t, \vec{z}), \qquad i=1,2,\hdots,q,\quad \vec{z}\in Z \subset \RR^q, \, t\in I\subset \RR,
	\end{align*}
	defined on the differential manifold $\dMan{M}_f$ with global coordinates  $(t,\vec{z})$ and the Cartan field $\vecf{\partial}_t = \tvec{t} + f^i(t, \vec{z}) \tvec{z^i}$. 	
	An automorphism $g:\dMan{M}_f \rightarrow \dMan{M}_f$ given by 
	\begin{align*}
		g: (t, \vec{z}) \mapsto (\tilde{t}, \tilde{z}) = (\theta(t, \vec{z}), \, \zeta(t,\vec{z})), \quad (t,\vec{z}),\, (\tilde{t},\tilde{\vec{z}})  \in I\times Z,
  \end{align*}
  is a Lie-B\"{a}cklund mapping iff the smooth functions  $\theta$ and $\zeta$ fulfill the defining equations \eqref{eq:LieBaecklund_definierendeGln}:
  \begin{align*}
  	\vecf{\partial}_t \zeta^i  - f^i\left(\tilde{t}, \tilde{z}\right) \vecf{\partial}_t \theta =0, \quad i=1,2,\hdots, q.
  \end{align*}
  This means nothing else but 
  \begin{align*}
  	f^i\left(\tilde{t}, \tilde{\vec{z}}\right) &= \frac{\pd{\zeta^i}{t}+ \pd{\zeta^i}{z^j} \, f^j(t,\vec{z})}{\pd{\theta}{t}+ \pd{\theta}{z^j} \, f^j(t,\vec{z})} = \frac{d\tilde{z}^i}{d\tilde{t}}, \quad i=1,2,\hdots,q,
  \end{align*}
  i.e.\ the differential equation is invariant w.r.t.\ any Lie-B\"{a}cklund mapping on $\dMan{M}_f$.  
\end{Bsp}

The definitions given above lack constructiveness in the sense that finding all symmetries means to solve the impossible problem of finding all Lie-B\"{a}cklund mappings of the form~\eqref{eq:LieBaecklundAbbildung_Form}  that preserve the Cartan distribution for a given system of differential equations, or equivalently, all contact transformations on $J^1 \pi_M$ that preserve the submanifold $\mathcal{S}$. 
However, it is well known that the situation eases if more structure within the family of symmetry transformations is assumed, namely to consider only transformations that are elements of a Lie group\footnote{Note that the additional structure narrows the considered class of symmetries down, e.g.\ discrete symmetries like reflection are excluded.} acting smoothly on $\mc{S}$. In order to follow this approach the next section recalls  some necessary facts about Lie groups. For details on Lie groups  and their application to differential equations see \cite{War83, Olv93}. 

\subsection{Lie symmetry groups and invariants}

		An {$r$-parameter Lie group} is a group $G$ which also carries the structure of an $r$-dimensional smooth manifold in such a way that both the group operation $m:G\times G \rightarrow G$, $m(g,h) = g\cdot h$, $g,h\in G$, and its inverse map $i:G\rightarrow G$, $i(g)=g^{-1}$, $g\in G$, are smooth maps.
Since the group operation carries a smooth inverse it defines a diffeomorphism on the manifold it acts on. From now on, only local groups acting on a smooth fibered $(1+q)$-dimensional manifold $\left(M, \pi, \mc{B}\right)$, $M\simeq\left(t, x^1, \hdots, x^q\right)=\vec{z}$, $\mc{B}\simeq(t)$, are considered. 

\begin{definition}
		Let $G$ be a Lie group with identity $e$ and $W\subset M$ an open set. A local \emph{transformation  group} $(\phi_g)_{g\in G}$ on $W$ is a smooth map $ (g,\vec{z})\in G \times W \mapsto \phi_g(\vec{z}) \in W $ such that $\phi_e(\vec{z})=\vec{z}$, $\vec{z}\in W$, and $\phi_{g_2}(\phi_{g_1}(\vec{z})) = \phi_{g_2\cdot  g_1}(\vec{z})$,  $g_1,g_2\in G$, $\vec{z}\in W$ where defined.	
\end{definition} 

In local coordinates one can write the action of the transformation group smoothly parametrized w.r.t.\ its $r$ group parameters $a^k$, $k=1,2,\hdots, r$,
\begin{align*}
	\phi_g: \RR^r \times W \rightarrow W,\, \vec{z}\mapsto \phi(\vec{z}; \vec{a}),\quad G\owns g\stackrel{\mathrm{loc}}{\simeq} \vec{a}=(a^1,a^2,\hdots, a^r)\in \RR^r.
\end{align*}
The action of all group elements on a given point defines the {orbit} of a transformation group, i.e.\ a minimal nonempty invariant subset $\mathcal{O}\subset W$. 
If the transformation group acts globally on $M$, the orbit $\mathcal{O}_z$ through $\vec{z}\in M$ is the set $\mathcal{O}_z=\{g\cdot \vec{z} | g\in G\}$. 
From the definition of an orbit it follows that a subset $W\subset M$ is $G$-invariant iff it is the union of orbits. Clearly, points belonging to the same orbit form an equivalence class w.r.t.\ their orbit membership. A  transformation group acts locally effectively if from $g \cdot \vec{z}=\vec{z}$ for all $\vec{z}\in W$ follows $g=e$. Further, a local transformation group acts {semi-regularly}  on $M$  if its orbit dimension is constant.  The group acts {regularly} if the action is semi-regular and for each point $ \vec{z} \in  W$ there exist arbitrarily small neighbourhoods $U_{\vec{z}}$ with the property that each orbit thorugh $\vec{z}$ of $G$ intersects $U_{\vec{z}}$ in a pathwise connected subset\footnote{As a necessary condition each orbit of $G$ is a
regular submanifold of $M$.}~\citep{Olv93}. The group acts free on $M$ if  for any two distinct $g, h \in G$ and $\forall \vec{z} \in W$, $g\cdot \vec{z} \neq h\cdot \vec{z}$, i.e.\ only the identity element $e$ has fixed points.

The transformation $\phi_g=(\phi^1,\phi^2,\hdots, \phi^{1+q})$ is the flow of its \emph{infinitesimal generators}
\begin{align} \label{eq:defining_equations_1}
	\vecf{v}_{k} &= \sum_{l=1}^{1+q} \, \left.\frac{d}{da^k} \phi^l(\vec{z}; \vec{a})\right|_{\vec{a}=0}\hspace{-2ex}\tvec{z^l}, \quad k=1,2,\hdots, r,
\end{align}
which are smooth vector fields spanning the tangent space of the orbits at the identity element ($\vec{a}=0)$. The key observation regarding transformation groups is that they are completely described by their associated Lie algebra, which in the case of effective group actions is isomorphic to the Lie algebra spanned by its {infinitesimal generators}. Consequently, instead of working with the local action $\phi_g$ it suffices to consider the $r$ infinitesimal generators\footnote{If the group action is not effective, the $r$ vector fields span a Lie algebra that is isomorphic to the Lie algebra of the effectively acting quotient group $G\backslash G_W$, $G_W=\cap_{\vec{z}\in W} G_{\vec{z}} =\left\{g\, | \, g\cdot \vec{z} =\vec{z},\, \forall \vec{z}\in W\right\}$ being the global isotropy subgroup of $G$. However, as the quotient group carries the same transformation properties, one can consider the quotient group in case that $G$ does not act effectively on $W$ \citep{Olv95}. }.

\begin{definition}
	A local transformation group $G$ acting on $\mathcal{S}$ is a symmetry group of the differential equation iff $g \cdot \mathcal{S}\subset \mathcal{S}$ for all $g\in G$.
\end{definition}

Now, let $W \subset M$ with adapted coordinates $(t,\vec{z})$. For each element $g\in G$ with $g\cdot W\subset W$ there exists an induced group action on the first jet bundle $J^1\pi$ corresponding to the first prolongation  $\prol{g}{1}$. As before in the general case Theorem~\ref{thm:baecklund} holds, and thus for more than one dependent variable any contact transformation on the jet manifold is derived as prolongation of a point symmetry. 
The symmetry condition $g \cdot \mathcal{S}\subset \mathcal{S}$ is fulfilled iff the group transformations only shifts solutions on $\mathcal{S}$, i.e.\ leaves $\mathcal{S}$ invariant. This is the case if the infinitesimal generators of the transformation group are tangent to $\mathcal{S}$. Combining B\"{a}cklund's theorem and Proposition~\ref{beh:tangentielles_vektorfeld} one arrives at the following symmetry conditions for a system in state representation and infinitesimal generators defined on $M$:
\begin{align} \label{eq:symmetry_defining_equations}
		\prol{\vecf{v}_k}{1}\left( F(\vec{p})\right) =0, \quad \forall \vec{p}\in\mc{S} \quad\text{and }\quad k=1,2,\hdots,r.
\end{align}
It turns out that the equations~\eqref{eq:symmetry_defining_equations} are necessary and sufficient for $\phi_g$ being a symmetry of $\mathcal{S}$. For a proof the reader is referred to \cite{Olv93}.

Recalling the idea of using invariant tracking errors for the control design, invariants of the considered symmetry group are of interest.

\begin{definition} An invariant is a real-valued function $I:M\rightarrow \RR$ for which $I(g\cdot \vec{z})=I(\vec{z})$, $\forall g\in G$ and $\forall \vec{z}\in M$, holds. Further, let $\prol{G}{k}$ be the $k$-th prolongation of $G$ acting on $J^k \pi$, $k\geq 1$. A differential invariant is a real-valued function $I:J^k \pi \rightarrow \RR$ for which $I(\prol{g}{k}\cdot \vec{p})= I(\vec{p})$, $\forall g\in G$ and $\vec{p}\in J^k\pi$, holds.
\end{definition}
From this definition it follows that every invariant function is constant along the orbits of $G$. The following theorem guarantees the existence of functionally independent invariants.

\begin{theorem}[\citealp{Olv95}]
	 Let $G$ be a Lie group acting semi-regularly on the $(q+1)$-dimensional manifold $M$ with $s$-dimensional orbits. Then, in a neighborhood $W$ of each point $\vec{z}_0\in W\subset M$, there exist $(1+q-s)$ functionally independent invariant functions $I_1(\vec{z}),\hdots,I_{1+q-s}(\vec{z})$.
\end{theorem}
The existence of invariant functions follows from the Frobenius theorem applied to a $s$-dimensional involutive subalgebra of the Lie algebra (compare \cite{War83}, Theorem~1.60). For a complete proof see also \cite{Ovs82}, §17.3. 
 Further, any other invariant (local) function $J$ can be expressed as function of $I_1(\vec{z}),\hdots,I_{1+q-s}(\vec{z})$, i.e.\ $J = H(I_1,\hdots,I_{1+q-s})$ for some suitable function $H$. 
Another important fact is that for any group $G$ acting locally effectively on $M$ with $r>\dim M$ there exists a number $\delta> 0$ such that the action of $\prol{G}{\delta}$ (the $\delta$-th prolongation of $G$) is locally free on $J^{\delta} M$ with orbit dimension $s=r$, i.e.\ the dimension of the orbits on the $\delta$-th jet manifold of $M$ is equal to the group dimension \citep[][§24.1]{Ovs82}.

\subsection{Computation of invariants by normalization}

A complete set of $(1+q-r)$ functionally independent invariants for $\prol{G}{\delta}$ can be constructed by following a normalization procedure first proposed in \cite{Kil89}, for a modern treatment see \cite{Olv99}:
\begin{enumerate}
	\item Find $\delta\geq 0$, such that $\Rg \left[ \pd{\prol{\phi_g}{\delta}\!\left(\vec{p};\vec{a}\right)}{\vec{a}}\right]_{a=0}\!\!\!\!\!\!\!=r$, for $\vec{p}\in W\subset J^{\delta}\pi$,  i.e.\ the prolonged action is locally free on $J^{\delta} \pi$ with orbit dimension $s=r$. 
	\item Select $r$ components $\left(\phi^1, \hdots, \phi^r\right)$  of $\prol{\phi_g}{\delta}$ such that the normalization equations
	\begin{align} \label{eq:normalisierung}
	 c_1 = \phi^{1}\!\left(\vec{p};\vec{a}\right),&& c_2 = \phi^{2}\!\left(\vec{p};\vec{a}\right), &&\hdots, &&c_r = \phi^{r}\!\left(\vec{p};\vec{a}\right), &&c_{\bullet} \in \RR, 
  \end{align}	 
		are regular w.r.t.\ the group parameter $\vec{a}$. The local solution yields an equivariant map $\vec{a}=\gamma\left(\vec{p}\right)$ called as \emph{moving frame}.
	\item Finally, a complete set of functionally invariants is obtained by using the map $\gamma$ in the remaining $(1+q-r)$ transformation equations,
		\begin{align*}
				 I_j=\left. \phi^{j+r}(\vec{p}) \right|_{a=\gamma\!\left(\vec{p}\right)}, \qquad j=1, \hdots, 1+q-r.
		\end{align*}
\end{enumerate}

\begin{figure}
	\centering
	\psfrag{a}{\hspace{-.5ex} $\vec{u}$}\psfrag{x}{\hspace{-.3ex}$\vec{v}$}   \psfrag{b}{\hspace{-.75ex} $\vec{w}$}
	\psfrag{Oa}{\hspace{-.5ex} $\mathcal{O}_{\vec{u}}$}\psfrag{Ox}{\hspace{-.5ex} $\mathcal{O}_{\vec{v}}$} \psfrag{Ob}{\hspace{-.5ex} $\mathcal{O}_{\vec{w}}$}
	\psfrag{Ux}{ $W$} \psfrag{M}{  $M$} 
	\psfrag{varphi}{$\psi$} \psfrag{varphi2}{ $\psi^{-1}$}
	\psfrag{nu}{ $z^{\Rmnum{1}}$} \psfrag{nu0}{\begin{parbox}{2ex}{\vspace{1ex}\hspace{-1ex} $\vec{z}_0^{\Rmnum{1}}$}\end{parbox} }
	\psfrag{xi}{ ${z}^{\Rmnum{2}} $} \psfrag{gam}{ $\gamma$}
	\includegraphics[width=.7\linewidth]{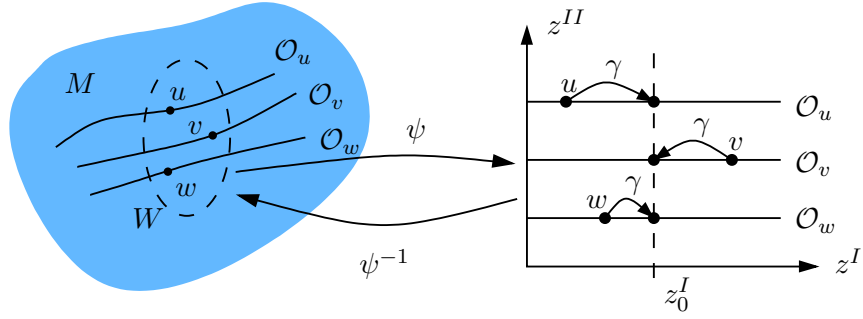}	
	\caption{Local foliation of the manifold $M$ generated by the group orbits, choice of a canonical element for each orbit, and meaning of $\gamma$ in terms of the group action in a local rectifying chart $(W, \psi)$}
	\label{fig:Blaetterung}
\end{figure}

The construction of the invariant functions becomes clear in local Euclidean coordinates for $J^{\delta} \pi$. Since the Lie algebra  (possibly of the prolonged group action) is involutive,  invocation of the Frobenius~Theorem allows one to choose local rectifying coordinates, as depicted in Figure~\ref{fig:Blaetterung}. For simplicity, assume that the group acts locally freely for $k=0$ with orbit dimension $r$ on $M$. In a suitable coordinate chart the first $r$ components $\vec{z}^{\Rmnum{1}}$ of $\vec{z}=\left(\vec{z}^{\Rmnum{1}}, \vec{z}^{\Rmnum{2}}\right)$ are the orbit coordinates whereas the remaining coordinates $\vec{z}^{\Rmnum{2}}$ are constant along the orbits. By solving the normalization equations~\eqref{eq:normalisierung} a canonical element is chosen for each orbit, $\vec{z}^{\Rmnum{1}}=\vec{z}^{\Rmnum{1}}_0$. For each point in $\vec{p}\in M$ the resulting map $\gamma: M \rightarrow G\stackrel{\mathrm{loc}}{\simeq}\RR^r$ defines a group element such that $\vec{p}$ is mapped to the associated canonical orbit element, i.e.\ the first $r$ coordinates coincide with $\vec{z}^{\Rmnum{1}}_0$. Consequently, the normalization introduces coordinates along the group orbits  w.r.t.\ the canonical element given by $\gamma$ -- it defines a \emph{rep\`{e}re mobile} (moving frame). As the normalization equations are equivariant w.r.t.\ the group action, their solution $\gamma$ enjoys the same property,  i.e.\ $\gamma(g\vec{z}) \cdot (g  \vec{z}) = \vec{z}_0 = \gamma\!\left(\vec{z}\right) \cdot \vec{z} \Rightarrow \gamma(g\vec{z})=\gamma(\vec{z}) g^{-1}$. Consequently, by using $\gamma$ in the remaining transformation equations for $\vec{z}$ in the last step of the normalization procedure, one obtains the (local) canonical element of the orbit, which is an invariant along the group orbits by construction.

\subsection{Computation of classical symmetries of nonlinear control systems} 

Together with a general ansatz for an infinitesimal symmetry of \eqref{eq:zustandssystem2} a set of linear first order partial differential equations can be derived from the symmetry conditions~\eqref{eq:symmetry_defining_equations} in order to compute the admitted symmetry group of a given system. Taking Theorem~\ref{thm:baecklund} into account, one arrives at the following ansatz for an infinitesimal generator on $M$
\begin{align} \label{eq:ansatzvektorfeld}
	\vecf{v} &= \xi(t,\vec{x}, \vec{u}) \tvec{t} +  \eta^i(t, \vec{x}, \vec{u})\tvec{x^i} +  \varphi^j(t, \vec{x}, \vec{u}) \tvec{u^j}.
\end{align}
Its first prolongation is given by (see \cite{Olv93} for a proof)\footnote{The $\prol{\vecf{v}}{1}$-notation is used instead of $j^1\vecf{v}$ since, in contrast to the prolongations before, the resulting transformation is in general not a bundle-morphism, i.e.\ it is not fiber-preserving, which would be the case for $\xi=\xi(t)$.}
\begin{align} \label{eq:symmetry_pr1}
\begin{split}
	\prol{\vecf{v}}{1} &= \vecf{v} +  \zeta^i\!\left(t, \vec{x}^{[1]}, 	\vec{u}^{[1]}\right) \tvec{\dot{x}^i} +  \psi^j\!\left(t, \vec{x}^{[1]}, \vec{u}^{[1]}\right) \tvec{\dot{u}^j}, \\
 & \text{with } \zeta^i = D_t \eta^i - \dot{x}^i D_t \xi, \quad \psi^j = D_t \varphi^j - \dot{u}^j D_t \xi,
\end{split} 
\end{align}
where $D_t$ denotes the total derivative w.r.t.\ $t$. Note that the bracket notation $[k]$ is used to denote time derivatives up to the $k$-th order. Defining the vectors $\vec{\eta}=\left(\eta^1,\eta^2,\hdots,\eta^n\right)^T$, $\vec{\zeta}=\left(\zeta^1,\zeta^2,\hdots,\zeta^n\right)^T$, $\vec{\varphi}=\left(\varphi^1,\varphi^2,\hdots,\varphi^m\right)^T$ and applying \eqref{eq:symmetry_pr1} to  \eqref{eq:submanifold_defining_equations} in combination with condition \eqref{eq:symmetry_defining_equations} yields
\begin{align*}
	\prol{\vecf{v}}{1}\left( F \right) &= D_t \vec{\eta} - \dot{\vec{x}} D_t \xi - f_{\vec{x}}  \vec{\eta} - f_{\vec{u}}  \vec{\varphi} = \vec{\eta}_t + \vec{\eta}_{\vec{x}} f + \vec{\eta}_{\vec{u}} \vec{\dot{u}} - f \left( \xi_t + \xi_{\vec{x}} f + \xi_{\vec{u}} \dot{\vec{u}} \right) - f_{\vec{x}} \vec{\eta} - f_{\vec{u}} \vec{\varphi}=0.
\end{align*}
Due to its linearity in $\dot{\vec{u}}$ this system splits into two parts:
\begin{align}  \label{eq:defining_equations_splitted}
	\begin{split}
	\vec{\eta}_{\vec{u}} - f \xi_{\vec{u}} &= 0, \\
	\vec{\eta}_t + \vec{\eta}_{\vec{x}} f  - f \left( \xi_t + \xi_{\vec{x}} f  \right) - f_{\vec{x}} \vec{\eta} - f_{\vec{u}} \vec{\varphi}&=0.
	\end{split}
\end{align}
In general, solving this system of linear partial differential equations is as difficult as solving the original differential equation. Hence, finding symmetries of a given nonlinear control system usually involves consideration of different ansatz vector fields reflecting some kind of intuition e.g.\ in terms of the dependence on certain system variables etc.

\begin{remark}
	Since the vector field~\eqref{eq:ansatzvektorfeld} generates a point transformation on $X\times U$ its prolongation to $J^1 \pi$ is contact by construction. Hence, all transformations generated by such a vector field naturally form Lie-Bäcklund maps of the contact distribution on $J^1 \pi$ formally qualifying them as potential symmetries.
\end{remark}

\section[Invariant control design]{Invariant control design for control systems with Lie symmetry}
\label{sec:invariant_control_design}

Consider a control system in state representation~\eqref{eq:zustandssystem2} and assume that it admits a local transformation group $\left(\varphi_g\times\psi_g\right)_{g\in G}$
\begin{align*}
		X \times U \owns \left(\tilde{\vec{x}}, \tilde{\vec{u}}\right) &= \left(\varphi_g(\vec{x}), \psi_g(\vec{x}, \vec{u})\right),\quad g\in G,
\end{align*}
where $\varphi_g$ is a local diffeomorphism and $\psi_g$ is regular w.r.t.\ $\vec{u}$ for all $\vec{x}\in X$ (coordinate change and regular state feedback), and the equality\footnote{Recalling Example~\ref{bsp:LieBaecklundForminvarianz} one observes that this is equivalent to $\left(\varphi_g\times\psi_g\right)_{g\in G}$ defining Lie-Bäcklund mappings on $X\times U$.}
\begin{align} \label{eq:equivariance}
			\pd{\varphi_g}{\vec{x}}\left(\vec{x}) f(\vec{x}, \vec{u}\right) &= f\left(\varphi_g(\vec{x}), \psi_g(\vec{x},\vec{u})\right)
\end{align}
holds for all $g\in G$ and $(\vec{x}, \vec{u})\in X\times U$. Note that $\varphi_g$ necessarily depends only on $\vec{x}$ since the time is not transformed (compare equation~\eqref{eq:defining_equations_splitted}). 
Further, assume that the known symmetry should be preserved under a feedback law that is to be designed in order to stabilize the output 
\begin{align*}
		\vec{y}=h(\vec{x}), \qquad h(\vec{x})=\left(h^1(\vec{x}), h^2(\vec{x}), \hdots, h^m(\vec{x})\right), \, h\in {C}^{\infty}(X),
\end{align*}
along some suitably planned smooth reference trajectory\footnote{For notational ease $\ol{\vec{y}}_d$ denotes the reference trajectory including time derivatives up to the required order.} $t \mapsto \left(\vec{y}_d(t), \dot{\vec{y}}_d(t), \hdots \right)=:\ol{\vec{y}}_d$. 

To this end, assume that the output $\vec{y}=h(\vec{x})$ has a well-defined vector relative degree $\vec{r}=\left(r^1, \hdots, r^m\right)$, $\sum_i r^i = n$, such that the decoupling matrix 
\begin{align*}
				\pd{}{\vec{u}} \begin{pmatrix} L_f^{r^1} y^1  & L_f^{r^2} y^2 & \cdots & L_f^{r^m} y^m \end{pmatrix}^T
\end{align*}
is generically regular on $X\times U$ \citep[for more details see e.g.\!][]{NvdS90, Isi95}. The smooth map $\zeta: X \rightarrow Z\subset \RR^m$, defining new coordinates $\vec{z}=\left(z^1, \hdots, z^m\right)$,
\begin{align*}
		&\zeta:  \quad z^1 = h^1(\vec{x}), \quad z^2 = h^2(\vec{x}), \quad \hdots, \quad z^{m} = h^m(\vec{x}),
\end{align*}
induces a local diffeomorphism 
\begin{align*}
	\Phi: \quad z^i_j = L_f^j h^i(\vec{x}),\quad j=0,\hdots, r^i-1, \quad i=1,\hdots, m,
\end{align*}
and one obtains the following system representation in $(\vec{z}, \vec{u})$-coordinates:
\begin{align} \label{eq:normalform}
	\begin{split}
	\begin{array}{rcl}
		\dot{z}^i_j &=& z^i_{j+1},\\
		\dot{z}^i_{r^i-1} &=& F^i(\vec{z}, \vec{u}),
	\end{array} \qquad \text{for } j=1,\hdots, r^i-2,\quad  i=1,\hdots, m.
	\end{split}
\end{align}
Since $\Phi$ defines a diffeomorphism on $X$ its first prolongation is a Lie-Bäcklund mapping w.r.t.\ the contact distribution of $J^1 \pi$. Therefore, the local symmetry group $(\varphi_g \times \psi_g)_{g\in G}$ induces the local symmetry group $\left(\ol{\varphi}_g\times\psi_g\right)_{g\in G}$ in $(\vec{z},\vec{u})$-coordinates with $\ol{\varphi}_g = \Phi \circ \varphi_g \circ \Phi^{-1}$. Defining the usual tracking errors $e^i=y^i-y^i_d$ w.r.t.\ the components of the input and the reference trajectory and imposing a linear time-invariant error dynamics
\begin{align} \label{eq:fehlerdynamik}
		e^i_{r^i} + c_i^{r^i-1} e^i_{r^i-1} + \cdots + c_i^{1} \dot{e}^i + c_i^0 e^i &= 0 = E^i\!\left(e^i,e^i_1, \hdots, \vec{e}^i_{r^i}\right),  \quad i=1,\hdots, m,
\end{align}
where $e^i_j$ denotes the $j$-th time derivative of the $i$-th error component, and using the model equations~\eqref{eq:normalform} one obtains a feedback law $\vec{u}=U\left(\vec{x}, \ol{\vec{y}}_d\right)$ realizing the desired error dynamics. However, rewriting the error dynamics in transformed $(\tilde{\vec{x}}, \vec{\tilde{u}})$-coordinates only results in the identical tracking behavior, iff each $E^i$, $i=1, \hdots, m$, is equivariant w.r.t.\ the local symmetry group, i.e.\ 
\begin{align*}
		j^{r^i} \ol{\rho}_g \circ E^i(e^i_{[r^i]}) = E^i\!\left(j^{r^i} \ol{\rho}_g \circ e^i_{[r^i]}\right) \qquad \text{for all } g\in G, \,i=1,\hdots, m,
\end{align*}
where $j^{r^i} \ol{\rho}_g \circ E^i$ and $j^{r^i} \ol{\rho}_g \circ e^i_{[r^i]}$  denote the induced action  on the error dynamics and the error, respectively. 
If these equalities hold, the stabilizing feedback law $U$ renders the system~\eqref{eq:normalform} invariant w.r.t.\ the induced symmetry group. Moreover, since symmetry is independent of the specific system representation, the symmetry of the original state representation is preserved under this feedback as well.

Equivariance is clearly achieved if the error is invariant w.r.t.\ the group action motivating an invariant feedback approach based on \emph{invariant tracking errors} as proposed in \cite{RR99, Rud03, MRR04}. The given feedback design assumes an induced group action on the output $\vec{y}$, which due to the special choice of the coordinates was given by $\ol{\varphi}_g$. In a more general setting this  motivates the following definition.

\begin{definition}[\citealp{MRR04}] An output $y = h(\vec{x}, \vec{u})$ is called $G$-compatible if there exists an induced transformation group $(\rho_g)_{g\in G}$ on $Y$ such that $h(g\cdot (\vec{x}, \vec{u})) = \rho_g \cdot h$ for all $g\in G$.
\end{definition}

The special distinction of $G$-compatible outputs becomes necessary as in the following the invariant tracking error will be defined in terms of the output $\vec{y}$ only. Consequently, the group action has to be restricted to $Y$ inducing the necessity of the co-distribution generated by the elements of the output map $h=\left(h^1,\hdots, h^m\right)^T$
\begin{align*}
				\Omega &= \Span\left\{ dh^1(\vec{x}), dh^2(\vec{x}), \hdots, dh^m(\vec{x}) \right \} \subset \CTRaum{X}{},
\end{align*}
to be invariant along the group orbits, i.e.\ $\Omega$ is closed w.r.t.\ the group action, 
\begin{align*} 
				\omega \in \Omega \quad \Rightarrow \quad L_{\vecf{v}_k} \omega \in \Omega \quad \text{for } k=1,2,\hdots, r\,,
\end{align*} 
where $L_{\vecf{v}_k}$ denotes the Lie derivative along the $k$-th infinitesimal generator $\vecf{v}_k$ of the group. This observation might be exploited if $G$-compatible outputs are not known as shown in  the following example.

\clearpage
\begin{Bsp}[Computation of $G$-compatible outputs]\label{Bsp:Fzg}

\begin{wrapfigure}{r}{.3\linewidth}
	\small 
	\psfrag{y}{$\vec{y}$}
	\psfrag{theta}{$\theta$}
	\psfrag{varphi}{$\varphi$}
	\psfrag{v}{$v$}
	\psfrag{l}{$l$}
	\psfrag{y1}{$z^1$}  	\psfrag{y2}{$z^2$}
	\includegraphics[width=\linewidth]{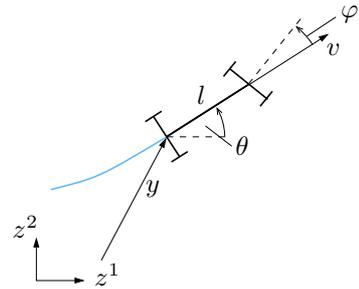}
	\caption{Kinematic car}
	\label{fig:Fzg}
\end{wrapfigure}

	Consider the well-known equations of the kinematic car (unicycle) 
\begin{align}\label{eq:kinematisches_Fzg_Modellgleichungen}
		\dot{\vec{z}} &= v \begin{pmatrix} \cos \theta \\ \sin \theta \end{pmatrix},\quad \dot{\theta} = \frac{v}{l} \tan \varphi ,
\end{align}	
with $\vec{z}=\left(z^1,z^2\right)\in \RR^2$ being the position of the rear axle center, $\theta$ describing the orientation of the car, and the inputs $v$ (velocity) and $\varphi$ (steering angle) (Figure~\ref{fig:Fzg}). The model equations are invariant w.r.t.\ the special Euclidean group SE(2) consisting of rotations and translations in the plane with infinitesimal generators
\begin{align*}
		\vecf{v}_1 &= -z^2 \tvec{z^1} + z^1 \tvec{z^2} + \tvec{\theta}, && \vecf{v}_2 = \tvec{z^1}, && \vecf{v}_3 = \tvec{z^2}.
\end{align*}
With $\vec{x}=\left(z^1, z^2, \theta\right)$ and $\vec{u}=\left(v, \varphi\right)$ the model equations are in state representation. Consider a two-dimensional smooth output $\vec{y}=\left(h^1(\vec{x}), h^2(\vec{x})\right)$, spanning the smooth co-distribution $\Omega=\Span \left\{ dh^1(\vec{x}), dh^2(\vec{x}) \right\}$. If $\Omega$ is invariant w.r.t.\ the group action this is necessarily also true for its annihilator $\Omega^{\perp}=\Span\left\{ \vecf{w} \right\}$ spanned by a smooth vector field $\vecf{w}(\vec{x})= w^1(\vec{x}) \tvec{x^1} + w^2(\vec{x}) \tvec{x^2} + w^3(\vec{x}) \tvec{x^3}$. The functions $w^i$, $i=1,2,3$, derive from $h$ using the conditions $\sprod{dh^i}{\vecf{w}}=0$.
Computing the Lie derivatives along the infinitesimal generators and re-using the annihilating property one obtains a set of six linear partial differential equations for $h^1$ and $h^2$
\begin{align*}
	\sprod{dh^i}{\liekl{\vecf{w}}{\vecf{v}_k}} = 0, \qquad i=1,2, \quad k=1,\hdots, 3.
\end{align*}
Using different ansatz functions for $h^1$, $h^2$ one can for instance find that any given pair of smooth function $h^i(z^1, z^2)$, $i=1,2$, forms a $G$-compatible output, showing that  $\vec{y}=(z^1, z^2)$ is a valid choice. Further, setting $h^i=h^i(z^1, \theta)$ yields the solutions
\begin{align*}
	h^1&= h^1(z^1,\theta), & h^2&= \mathrm{const.}; && &h^1 &= h^1(z^1),        &h^2 &= h^2(z^1); \\
  h^1&= h^1(\theta),     & h^2&= h^2(\theta);     && &h^2 &= h^2(z^1,\theta), &h^1 &= h^1(-h^2(z^1,\theta)), 
\end{align*}
 which are not suitable for control purposes.
\end{Bsp}

\begin{remark}[$G$-compatible flat outputs] 
	In \cite{MRR97} the notion of \linebreak \emph{symmetry-preserving} flat outputs is given which are $G$-compatible flat outputs, i.e.\ the considered symmetry group induces a local group of transformations acting on the flat output. Further, the flat output $\ol{\vec{y}}$ and its image  $\tilde{\ol{\vec{y}}}$ under SO(2) 
\begin{align*}
	\ol{\vec{y}}&= \left( z^1, z^2 + \dot{z}^1\right), \qquad \tilde{\ol{\vec{y}}} = \left(z^1 \cos a - z^2 \sin a, z^1 \sin a + z^2\cos a + \dot{z}^1 \cos a - \dot{z}^2 \sin a \right)
\end{align*}
for the kinematic car are given as an example for a flat output  that does not preserve the symmetry w.r.t.\ SO(2) since $\tilde{\ol{\vec{y}}}$ is clearly not transformed on $Y$ only. However, since symmetry is invariant under Lie-Bäcklund equivalence, any symmetry transformation acting on a particular flat output induces symmetry transformations on any other flat output\footnote{In fact, any contact transformation on a flat output forms a symmetry, and hence induces a symmetry transformation on any other flat output of the system.}. Consequently, every flat output is symmetry-preserving in the suitable context. For the example above this can be seen as follows. As recalled in  Remark~\ref{rem:flatness} every flat system is Lie-Bäcklund equivalent to a trivial system, i.e.\ an infinite chain of integrators of width $m$. For the well-known flat output $\vec{y}=\left(z^1, z^2\right)$ one obtains the differential manifold $\dMan{M}_{{y}}=\left( \RR\times \RR_2^{\infty}, \Span\{\vecf{v}_y\}\right)$, with the Cartan vector field $\vecf{v}_y=\tvec{t} + y^1_1 \tvec{y^1} + y^2_1 \tvec{y^2} + y^1_2 \tvec{y^1_1} + y^2_2 \tvec{y^2_1} + \cdots $. This differential manifold is Lie-Bäcklund equivalent to the differential manifold $\dMan{M}_f=\left(\RR\times \RR^3 \times \RR_2^{\infty}, \Span\{\vecf{v}_f\}\right)$ defined by the model equations via the Cartan vector field $\vecf{v}_f = \tvec{t} +  v \cos\theta \tvec{z^1} +  v\sin\theta \tvec{z^2} + \frac{v}{l}\tan \varphi \tvec{\theta}  + \dot{v} \tvec{v} + \dot{\varphi} \tvec{\varphi} + \ddot{v} \tvec{\dot{v}} + \cdots$, which can be seen by expressing $z^1, z^2, \theta$ and $v,\varphi$ in terms of the flat output and its time derivatives,
\begin{align*}
	z^1 = y^1, && z^2 = y^2, && \theta = \arctan\frac{y^2_1}{y^1_1}, && v= \pm \sqrt{(y_1^1)^2+(y^2_1)^2}, && \varphi = \pm\arctan \left( l \frac{y^2_2 y^1_1 - y^1_2 y^2_1}{\left((y^1_1)^2+(y^2_1)^2 \right)^{\frac{3}{2}}}\right).
\end{align*} 
On the other hand, 	the differential manifolds $\dMan{M}_{y}$ is Lie-Bäcklund equivalent to $\dMan{M}_{\ol{y}}$ (which is defined analogously), i.e.\ both manifolds are transformed into each other using the mappings 
\begin{align*}
	& \Phi: \dMan{M}_y \rightarrow \dMan{M}_{\ol{y}}, \, \ol{y}^1 = y^1,\, \ol{y}^2 = y^2 + {y}^1_1,\hdots \quad \text{and}\quad  \Psi: \dMan{M}_{\ol{y}} \rightarrow \dMan{M}_{{y}}, \, {y}^1 = \ol{y}^1,\, {y}^2 = \ol{y}^2 - {\ol{y}}_1^1,\hdots\\
	& \pushforward{\Phi}(\vecf{v}_{{y}}) = \tvec{t} + y^1_1 \tvec{\ol{y}^1} + (y^2_1 + y^1_2) \tvec{\ol{y}^2}  + y_2^1 \tvec{\ol{y}^1_1} + (y_2^2 + y^1_3) \tvec{\ol{y}^2_1} + \cdots\\
																	 &= \tvec{t} + \ol{y}^1_1 \tvec{\ol{y}^1} + (\ol{y}^2_1 - \ol{y}^1_2 + \ol{y}^1_2) \tvec{\ol{y}^2} + \ol{y}^1_2 \tvec{\ol{y}^1_1} + (\ol{y}^2_2 - \ol{y}^1_3 + \ol{y}^1_3) \tvec{\ol{y}^2_1} + \cdots = \vecf{v}_{\ol{y}}.
\end{align*}
Hence, the symmetry transformation $\tilde{\vec{y}} = \left(y^1\cos a -y^2\sin a, y^1\sin a + y^2 \cos a\right)$, $a\in [0, 2\pi)$, induces the transformation
\begin{align*}
	\RR_2^{\infty} \rightarrow \RR_2^{\infty}: \qquad \tilde{\ol{y}}^1 = \ol{y}^1\cos a -(\ol{y}^2- \ol{y}_1^1) \sin a, && \tilde{\ol{y}}^2 = \ol{y}^1\sin a + \ol{y}^2\cos a - (\ol{y}^2_1 - \ol{y}^1_2)\sin a, && \hdots\, 
\end{align*}
which preserves the Cartan distribution, and therefore is a symmetry w.r.t.\ the system representation $\dMan{M}_{\ol{y}}$, showing, as claimed before, that $\ol{\vec{y}}$ is symmetry-preserving as well. However, as both system representations are connected via an infinite Lie-Bäcklund isomorphism (i.e.\ not a diffeomorphism/ point transformation), the symmetry is not classical but a generalized symmetry in the aforementioned sense.
\end{remark}

Based on the notion of $G$-compatible outputs invariant tracking errors can be readily defined.

\begin{definition}[\citealp{MRR04}] A smooth mapping $\left(\vec{y},\ol{\vec{y}}_d\right)\mapsto I\left(\vec{y},\ol{\vec{y}}_d\right)$ is a $G$-compatible tracking error if $\vec{y} \mapsto I\left(\vec{y},\ol{\vec{y}}_d\right)$ is locally invertible w.r.t.\ $\vec{y}$ along the reference trajectory, if $I\left(\vec{y}_d, \ol{\vec{y}}_d\right)=0$  holds for all $\vec{y}_d\in{Y}$, and if $I$ is an invariant w.r.t.\ the induced group action  $(\rho)_{g\in G}$ on $Y$.
\end{definition}

\begin{proposition}[\citealp{MRR04}] Let $G$ be a regular, locally effective Lie group and $\vec{y}$ a $G$-compatible output. Then there locally exists a $G$-compatible tracking error.
\end{proposition}
 Recalling the fact that every locally effective group can be prolonged up to some large enough order $\delta$ to become locally free with orbit dimension $r$ and invoking the normalization algorithm on the reference trajectory, a local moving frame $g=\gamma(\ol{\vec{y}}_d)$ is obtained. Using $\gamma$ in the remaining transformation equations for $\vec{y}$ yields a set of $m$ functionally independent invariants $I_j\left(\vec{y}, \ol{\vec{y}}_d\right)=\left.\rho_g \cdot y^j \right|_{g=\gamma(\ol{\vec{y}}_d)}$. Further, invariant tracking errors are given by $e^j = I_j\left(\vec{y}, \ol{\vec{y}}_d\right) - I_j\left(\vec{y}_d, \ol{\vec{y}}_d\right)$, $j=1,\hdots, m$. 

Once a complete set of invariant tracking errors has been derived, invariant feedback laws, i.e.\ feedbacks that preserve the symmetry, can be designed using well-known design approaches such as input-output linearization, integrator backstepping, or sliding mode control.

\begin{Bsp}[Invariant tracking errors obtained by normalization] 
	The following tracking errors can be derived by normalization:
\begin{itemize}
	\item \emph{kinematic car, planar motion of a rigid body:} The symmetry w.r.t.\ SE(2) is preserved when a moving frame attached to the reference trajectory for a flat output $\vec{y}$ is used, e.g.\ \linebreak $\vec{e}=\left(\sprod{\vec{y}-\vec{y}_d}{\vec{\tau}}, \sprod{\vec{y}-\vec{y}_d}{\vec{\nu}}\right)$, with $t \mapsto \vec{y}_d(t)$ denoting a smooth reference trajectory with well-defined tangent and normal vector $\vec{\tau}$ and $\vec{\nu}$, see for instance \cite{Woe98, MRR04, RF03}.
	
	\item \emph{attitude control of a rigid body (e.g.\ satellite):} Let $R,R_d \in \mathrm{SO}(3)$ denote two rotational matrices describing the current and desired orientation of a rigid body w.r.t.\ some inertial frame. Using the normalization equation $R_g R_d = R_e =I$, with $R_g$ being an element of SO(3) representing the group action and $I$ denoting the identity matrix, one obtains the map $R_g=R_{\gamma(R_d)}=R_d^T$ yielding the invariant error $R_d^T R$. Using a quaternion approach to parametrize the rotational matrices one arrives at the well-known invariant rotational error $\vec{e} = \vec{q}_d^{-1} \vec{q}$ \citep{WK91}. More generally, if one deals with a system on Lie groups this construction holds generally, i.e.\ the invariant error is given by $g_d^{-1} \cdot g$ derived from $g \cdot g_d = e$ with $e$ being the identity element of the Lie group \citep[see for instance][]{MBD06}.
	
	\item \emph{scaling group:} Using a relative error when a symmetry acts as a scaling  on the output yields an invariant error $\vec{e}=\frac{\vec{y}}{\vec{y}_d}-1$, see for instance \cite{RR99}.
\end{itemize}	
\end{Bsp}

\subsection{Structural consequences: system representation of reduced order}

As pointed out in \cite{GM83, GM85}, the existence of Lie symmetries imposes consequences on the underlying structure of the control system, which is closely related to the notion of a controlled invariant distribution \citep[compare][]{NvdS1982, NvdS1985, Isi95}. In order to give an example of further application of the normalization procedure it is necessary to recall some facts regarding the structural consequences of Lie symmetries. Consider a so-called \emph{state-symmetry}, i.e.\ a symmetry group acting locally freely only on the state variables with orbit dimension $r$. Since the infinitesimal generators form a Lie algebra their Lie brackets w.r.t.\ each other can be rewritten in terms of infinitesimal generators and the structure constants $c_{ij}^k\in \RR$ of the Lie group,
\begin{align*}
		\liekl{\vecf{v}_i}{\vecf{v}_j} = c_{ij}^1 \vecf{v}_1 + c_{ij}^2 \vecf{v}_2 + \cdots + c_{ij}^r \vecf{v}_r.
\end{align*}
Note that the Lie algebra forms an involutive distribution on $X$ and hence, there exist a local chart for $X$ such that the orbits of $G$ form submanifolds with $x^i=\mathrm{const.}$, $i=r+1, \hdots, n$, i.e.\ the group acts only along the first $r$ coordinate directions. The construction of such a chart is for instance given in \cite{War83} along with the proof of Theorem~1.60. In this new chart the infinitesimal generators have the form 
\begin{align*}
		\vecf{v}_i &= \varphi_i^1\!\left(x^1, \hdots, x^r\right) \tvec{x^1} + \varphi_i^2\!\left(x^1, \hdots, x^r\right) \tvec{x^2} + \cdots + \varphi_i^r\!\left(x^1, \hdots, x^r\right) \tvec{x^r}, \quad i=1,\hdots,r.
\end{align*}
Consequently, the symmetry conditions $\prol{\vecf{v}_i}{1}\!\left(\dot{\vec{x}}- f(\vec{x}, \vec{u})\right)$ on $\mc{S}$ read
\begin{align} \label{eq:Struktur_Lie-Klammer_Kommutator}
		&\left(\begin{array}{c:c}
										\begin{array}{ccc} 
										\pd{\varphi_i^1}{x^1} &  \cdots & \pd{\varphi_i^1}{x^r}  \\
											\vdots              &   \vdots &   \vdots  \\
									  \pd{\varphi_i^r}{x^1} &  \cdots & \pd{\varphi_i^r}{x^r} \\
									  \end{array}
									  & 0_{r\times(n-r)} \\ \\[-1.5ex] \hdashline \\[-1.5ex]
									   0_{(n-r)\times r} & 0_{(n-r)\times(n-r)} 
			\end{array}	\right)			
									  \begin{pmatrix} f^1 \\ f^2 \\ f^3 \\ \vdots \\ f^n \end{pmatrix} -  
			\left[ \pd{f}{\vec{x}}\right]
			\left(\begin{array}{c}						  									  
									   	\varphi_i^1 \\  \vdots \\ \varphi_i^r \\ \hdashline  0 \\ \vdots \\0 
		  \end{array} \right)  = \liekl{\vecf{v}_f}{\vecf{v}_i} = 0, 				
\end{align}
for $i=1,\hdots, r$, i.e.\ the infinitesimal generators commute with the vector field $\vecf{v}_f= f^i \tvec{x^i}$. From this it follows that the Lie algebra $\basis{g} = \Span\left\{ \vecf{v}_1, \hdots\, \vecf{v}_r \right\}$  is 
invariant along the flow of $\vecf{v}_f$, i.e.\ $\liekl{{\vecf{v}}}{\vecf{v}_f} \in \basis{g}$ for  all $\vecf{v} \in \basis{g}$ -- it forms a controlled invariant distribution of the control system. Examining the Lie brackets of the elements $\basis{g}\owns \vecf{e}_i = \tvec{x^i}$, $i=1,\hdots, r$,  with the vector fields $\vecf{v}_f$ yields the conditions
\begin{align*}
	\liekl{\vecf{v}_f}{\vecf{e}_i} =  -\begin{pmatrix}\pd{f^1}{x^i} & \pd{f^2}{x^i} & \cdots & \pd{f^n}{x^i} \end{pmatrix}^T = \lambda_i^1\!\left(x^1, \hdots, x^r\right)  \tvec{x^1} + \cdots + \lambda_i^r\!\left(x^1, \hdots, x^r\right) \tvec{x^r}
\end{align*}
which have to be fulfilled for suitable $\lambda_i^k$. Hence, the functions $f^{r+1}, \hdots, f^n$ have the form  $f^k=f^k\!\left(x^{r+1}, \hdots, x^n, u^1, \hdots, u^m\right)$, $k=r+1,\hdots, n$ in this particular coordinate chart, i.e.\ one arrives at the system representation  
\begin{align*} 
			\dot{\vec{\xi}} = F_1\!\left(\vec{\xi}, \vec{u}\right),  && \dot{\vec{\eta}}= F_2\!\left(\vec{\xi}, \vec{\eta}, \vec{u}\right), \quad \text{with  } \vec{\xi}=\left(x^{r+1}, \hdots, x^n \right), \,\vec{\eta}=\left(x^1, \hdots, x^r\right).
\end{align*}
From this local representation a possible decomposition of the control system into two parts can be observed:
\begin{itemize}
	\item a $(n-r)$-dimensional system $\dot{\vec{\xi}}=F_1(\vec{\xi}, \vec{u})$ describing the motion transverse to the orbits (i.e.\ ``from orbit to orbit"), and
	\item an $r$-dimensional system $\dot{\vec{\eta}}=F_2(\vec{\xi}, \vec{\eta}, \vec{u})$ describing the motion along the orbits\footnote{Since $\vec{\eta}$ locally defines coordinates on $G$, i.e.\ $g\stackrel{\mathrm{loc}}{\simeq} \vec{\eta}$, this part of the motion can be interpreted as motion on the Lie group. In fact, if the decomposition exists globally one can express this part as dynamical system on $G$ as mentioned in Remark~\ref{rem:global_decomposition}.}. 
\end{itemize} 
Depending on the control problem, it may be sufficient to consider the $(n-r)$-dimensional reduced order system, as the motion along the orbits is completely determined by the evolution of $\vec{\xi}$ (and possibly some remaining components of the input), see also \cite{ZZ92}. The coordinate transformation yielding the $(n-r)$ coordinates $\xi^1, \hdots, \xi^{n-r}$ can be constructed using the normalization procedure. Solving $r$ normalization equations (possibly renumbering the components of the state)
\begin{align*}
			\varphi_g\cdot x^1 &= c^1, && \varphi_g\cdot x^2 = c^2, && \hdots \,, && \varphi_g\cdot x^r = c^r, \quad c^i\in \RR, 			
\end{align*}
yields a moving frame $g=\gamma\!\left(x^1,\hdots, x^r\right)$. By the usual procedure the local coordinate transformation is obtained from setting
\begin{align*}
	\xi^1= \Xi^1(\vec{x}) = \left. \varphi_g \cdot x^{r+1}\right|_{g=\gamma\left(x^1,\hdots, x^r\right)}, && \hdots\,, && \xi^{n-r}= \Xi^{n-r}(\vec{x}) =  \left. \varphi_g \cdot x^{n}\right|_{g=\gamma\left(x^1,\hdots, x^r\right)},
\end{align*}
allowing to introduce the reduced order representation
\begin{align} \label{eq:Red:einsetzen}
	\dot{\xi}^i &= \left. \pd{\Xi^i}{\vec{x}}(\vec{x}) f(\vec{x}, \vec{u}) \right|_{\vec{x}=\left(c^1, \hdots, c^r, \Xi^{-1}(\vec{\xi})\right)},
\end{align}
with $\Xi=\left(\Xi^1, \hdots, \Xi^{n-r}\right)$. An application of this procedure is given in Example~\ref{Bsp:Red:Fzg}.

\begin{remark} \label{rem:global_decomposition}
	For a global reduction it is necessary to construct a global decomposition of $\vecf{v}_f$ along and transversal to the group orbits. A suitable geometric formulation can be given by using the bundle $(X\times U, \pi_G, X/G \times U)$,  $\pi_G: X \times U  \rightarrow X/G\times U $, $\pi_G(\vec{x}, \vec{u}) = \left(\mc{O}_{\vec{x}}, \vec{u}\right)$. Assuming that the group acts (globally) freely and properly on $X$ with orbit dimension $r$, the orbit space $X/G$ is in fact a smooth manifold and $\pi_G$ is a submersion, justifying the bundle construction (cf.~\cite{AM87}, Proposition 4.1.23). Therefore, a unique decomposition of  $\vecf{v}_f$ into a horizontal vector field being a section in $\TRaum{(X/G)}{}$ and a vertical vector field is possible based on the global foliation of $X$ generated by the group orbits. Hence, the quotient system on $X/G\times U$ can be defined globally. Further, the remaining $r$-dimensional systems can be rewritten as system on the Lie group $G$.  For details the reader is referred to \cite{GM85} and \cite{NvdS1982, NvdS1985}. 
\end{remark}

\begin{Bsp}[Kinematic car: reduced order state representation] 
\label{Bsp:Red:Fzg}

\begin{wrapfigure}[35]{r}{.35\linewidth}
	\centering	
	\vspace{-6ex}
	\small 
	\psfrag{y1}{$y^1$} \psfrag{y2}{$y^2$}
	\psfrag{z1}{$z^1$} \psfrag{z2}{$z^2$}
	\psfrag{thetadot}{$\dot{\theta}$}
	\psfrag{H}{$\vec{y}$}
	\psfrag{phi}{$\varphi$}
	\psfrag{theta}{$\theta$}
	\psfrag{v}{$v\tau $}
	\psfrag{l}{$l$}
	\includegraphics[width=\linewidth]{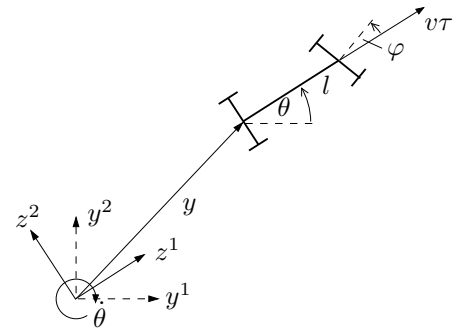}
	\caption{Interpretation of the $\vec{z}$-coordinate frame}
	\label{fig:Bsp:Red}
\end{wrapfigure} 

Based on the rotational invariance of the model equations~\eqref{eq:kinematisches_Fzg_Modellgleichungen} consider the normalization equation
\begin{align*}
			\tilde{\theta} =	\theta + a = 0 \quad \Rightarrow \quad a=\gamma(\theta)=-\theta.
\end{align*}
Solving for the group parameter $a$ (angle of rotation) and using the resulting map $\gamma$ in the transformation equations for $\vec{y}$ yields the  new (rotational invariant) coordinates
\begin{align*}
		\vec{z} &= R_{(-\theta)} \vec{y} = R_{\theta}^T \vec{y},  
\end{align*}
where $R_{\theta}\in \mathrm{SO}(2)$ denotes the rotational matrix w.r.t\ $\theta$. Equation~\eqref{eq:Red:einsetzen} yields the equations of motion in $\vec{z}$-coordinates:
\begin{align*}
	\dot{\vec{z}} &= \left[ \begin{pmatrix} \cos \theta & \sin \theta \\ -\sin\theta & \cos\theta \end{pmatrix} \begin{pmatrix} \cos\theta \\ \sin\theta \end{pmatrix}  v \right.\\ & \left. + \dot{\theta} \begin{pmatrix} -\sin\theta & \cos\theta \\ -\cos\theta & -\sin\theta  \end{pmatrix}  \begin{pmatrix} \cos\theta & -\sin\theta \\ \sin\theta & \cos \theta  \end{pmatrix}  \vec{z}  \right]_{\theta=0} 
	             = v \begin{pmatrix} 1  \\ 0 \end{pmatrix} + \frac{v}{l} \tan \varphi \begin{pmatrix} z^2   \\ -z^1 \end{pmatrix}.
\end{align*}
This can be interpreted as follows: The new $\vec{z}$ describe the position of the rear axle midpoint w.r.t.\ a coordinate frame that is rotated about the origin by $\theta$ and is rotating with angular velocity  $-\dot{\theta}$ ($\tilde{\theta} \equiv 0$). Since the translation along the tangent direction with velocity $v$ and the rotation of the coordinate frame are decoupled motions, the motion of $\vec{z}$ is described by the superposition of both (Figure~\ref{fig:Bsp:Red}). The normalization equation also delivers the equation for the motion along the group orbit. From $0 + a = \theta$ one observes that is is described by the third model equation, i.e.\ $g\simeq a$, 
$\dot{a} = \frac{v}{l}\tan\varphi$. \newline 
\end{Bsp}

\section[Application of invariant feedback design]{Application of invariant feedback design to a predator-prey bioreactor}
\label{sec:Beispiel}

Turning to the application of invariant tracking errors, consider the model of a continuously operated (chemostat) bioreactor containing two populations $P$ and $B$ of microorganisms of concentration $p$ and $b$ respectively. The predator population $P$ feeds upon the prey population $B$ whereas the population $B$ lives on some substrate $S$ of concentration $s$. The reactor is fed by a medium containing the substrate of concentration $s_F$ at the  dilution rate $D=\frac{Q}{V}$, where $Q$ and $V$ denote the flow rate and the volume of the reactor content.
A simplified model describing the evolution of the two populations reads \citep{Pav85}\\[-2ex]
\begin{align} \label{eq:bioreactor:alle}
	\begin{split}
	\dot{p} &= -Dp + \nu(b) p \\
	\dot{b} &= -Db + \mu(s) b - {\alpha} \nu(b) p \\
	\dot{s} &= D\left(s_F - s\right)  - {\beta} \mu(s) b
	\end{split}
\end{align}
where $\alpha$ and $\beta$ are the yield coefficients for production of both populations, and $\nu$ and $\mu$ are bounded, monotonically increasing functions modeling the specific growth rates of $P$ and $B$. In the following two types of growth models are considered (Figure~\ref{fig:Beispiel_Bioreaktor_Kinetiken}):

\begin{figure}
	\centering
	\psfrag{nu_max}{\hspace{1.5ex}\footnotesize $\nu_m$}
			\psfrag{nu_ext}{\hspace{1.5ex}\footnotesize $\overline{\nu}$}
			\psfrag{b_max}{\hspace{1ex}\begin{parbox}{2ex}{\vspace{1ex} \footnotesize $\overline{b}$}\end{parbox}}
			\psfrag{nu1}{\begin{parbox}{2ex}{\vspace*{6ex} \footnotesize ${\nu_{\,\text{M}}}$}\end{parbox}}
			\psfrag{nu2}{\begin{parbox}{2ex}{\vspace*{6ex} \footnotesize ${\nu_{\,\text{H}}}$}\end{parbox}}
			\psfrag{b}{\footnotesize $b$} 
	\includegraphics[width=.7\linewidth]{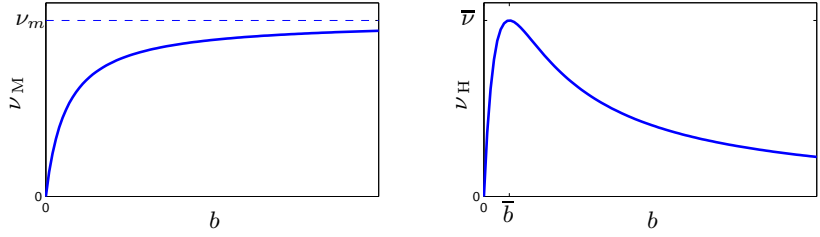} 
	\caption{Growth rate as function of the concentration $b$ for the Michaelis-Menten kinetic (left) and for the Haldane kinetic (right)}
	\label{fig:Beispiel_Bioreaktor_Kinetiken}
\end{figure}

\begin{enumerate}
\item Incorporating asymptotic saturation effects occurring at higher concentrations the well-known  Michaelis-Menten kinetic \citep{MM13} is applied, given by 
	\begin{align*}
				\nu_{\,\mathrm{M}}(s) = \nu_m \frac{b}{L+b},
	\end{align*}
with $\nu_m>0$, $L>0$ denoting the maximal growth rate and saturation constant, respectively.

\item If in addition to saturation also inhibition effects occur, the Haldane kinetic \citep{Hal30}  
\begin{align*}
				 \nu_{\,\mathrm{H}}(b) &=	\frac{\nu_{\text{m}} b}{b + K_S + K_I b^2 } 
\end{align*}
is used with maximal growth rate $\overline{\nu}= \nu_{\,\mathrm{H}}({\overline{b}}) =  \nu_m \frac{1}{1+2\sqrt{K_S K_I}}$ at $\overline{b}=\sqrt{\frac{K_S}{K_I}}$ parametrized using the saturation constant  $K_S$ and the inhibition coefficient $K_I$.
\end{enumerate}

\noindent The system~\eqref{eq:bioreactor:alle} possesses three kinds of equilibrium points:  washout of both populations for  $p=b=0$, $s=s_F$, washout of the predator population for  $p=0$, $b>0$, $s_F>s>0$, and coexistence of the two populations for $p>0$, $b>0$, $s_F>s>0$. The reactor is  operated at the third kind of equilibria which allows the two populations to coexist. 
For a detailed analysis of the equilibrium points, please see \cite{Pav85}.  

In the following, it is assumed that the control objective is to stabilize the reactor around equilibrium points in order to reject unmodeled disturbances using the output $\vec{y}=(p,b)$ and the input $\vec{u}=(D, s_F)$. Further, the growth rate for $b$ is modeled using the Michaelis-Menten kinetic, i.e.\ 
\begin{align*}
	\mu(s) &= \mu_{\mathrm{m}} \frac{s}{K+s}, \qquad K, \mu_{\mathrm{m}}>0.
\end{align*}
The stabilizing control should be realized independently of the growth kinetic $\nu$ used to model the growth of $P$. Consequently, the feedback design intends to render the feedback invariant w.r.t.\ the change from one kinetic to the other, i.e.\ smooth changes of the inhibition coefficient. To this end, the system equations~\eqref{eq:bioreactor:alle} are augmented by
\begin{align*}
	\dot{K_I} &= 0, && \nu(b) = \frac{\nu_{\text{m}} b}{b + K_S + K_I b^2 }.
\end{align*}
Changing from one growth kinetic to the other can then be interpreted as point transformation 
\begin{align*}
	\tilde{\nu}(b; a) =\frac{\nu_{\text{m}} b}{b + K_S + (K_I+a) b^2 }, \quad a \in [-K_I,0],
\end{align*}
acting on the inhibition coefficient $K_I$, where $a=0$ realizes the Haldane kinetic and $a=-K_I$ realizes the  Michaelis-Menten kinetic.
Using the invariance condition~\eqref{eq:equivariance} and choosing $p$ and $b$ to be invariant under the model change, i.e.\  $\tilde{p}=p$ and $\tilde{b}=b$, one obtains determining equations for the induced symmetry transformation on the substrate concentration $s$, the dilution rate $D$, and the feed concentration $s_F$:
\begin{align*}
		p\left({\nu_{\,\text{H}}}(b)-D\right) &= p \left(\tilde{\nu}(b;a) - \tilde{D}\right) && \Rightarrow		{\tilde{D}=\delta(b,D;a)}, \\
		b \left(\mu(s) - D\right) - \alpha {\nu_{\,\text{H}}}(b) p &=  b \left(\mu(\tilde{s}) - \tilde{D}\right) - \alpha \tilde{\nu}(b;a) p && \Rightarrow {\tilde{s}=\sigma(p,b,s;a)}, \\
		 \pd{\sigma}{p} p\left({\nu_{\,\text{H}}}(b)-D\right) + &\pd{\sigma}{b}\left(  b \left(\mu(s) - D\right) -  \alpha {\nu_{\,\text{H}}}(b) p \right)+  \\
		  \pd{\sigma}{s} \left( D\left(s_F - s\right)  - {\beta} \mu(s) b \right)  &= \tilde{D}\left(\tilde{s}_F - \tilde{s}\right)  - {\beta} \mu(\tilde{s}) \tilde{b} && \Rightarrow		{\tilde{s}_F=\sigma_F(p,b,s,s_F;a)}.
\end{align*}
Consequently, in order to render the mapping 
\begin{align*}
		\tilde{K}_I = K_I + a, && \tilde{p} = p, && \tilde{b}=b, && \tilde{s}=\sigma(p,b,s;a), 
\end{align*}
a symmetry of the reactor equations an induced feedback $(D, s_F)=(\delta(b,D;a), \sigma_F(p,b,s, s_F;a))$ has to be applied\footnote{In \cite{SB05} symmetries that are obtained by suitable feedback are denoted as \emph{controlled symmetries}.}. However, the system equations are locally feedback equivalent to a system for which the demanded symmetry acts only on the state variables $\vec{x}=\left(p,b,s\right)^T$. Again, this representation can be obtained by using the normalization algorithm. Normalizing the induced transformation $\tilde{s}=\sigma(p,b,s;a)=\mathrm{const}.$ yields a moving frame $a=\gamma(p,b,s)$, from which one derives the regular feedback (diffeomorphism w.r.t.\ $\vec{u}$),
\begin{align*}
	\vec{v} = V(\vec{x}, \vec{u}) = \left.\left(\delta(b,D;a), \sigma_F(p,b,s, s_F;a)\right)\right|_{a=\gamma(p,b,s)},
\end{align*}
with the new input $\vec{v}$ being invariant along the group orbits. 
As the symmetry has been chosen to leave the output $\vec{y}=(p,b)$ invariant, the usual set-point error $\vec{e}=\vec{y}- \vec{y}_0$ w.r.t.\ an equilibrium point $\vec{y}_0=\mathrm{const.}$ is a suitable choice for an invariant feedback design based on input-output linearization. 

\begin{remark} Following the spirit of the given example it is possible to define errors that are invariant w.r.t.\ set-point changes, i.e.\ stabilizing feedback laws that realize an identical error dynamics around a predefined set of set-points. For a set-point invariant integrator backstepping design please see \cite{CR10}. 
\end{remark}

\section{Conclusion}

 The advocated geometric approach to ordinary differential equations allows an intuitive interpretation of symmetries and the notion of Lie-Bäcklund equivalence.  As an important special subclass of classical symmetries Lie symmetries of nonlinear control systems have been considered, including their computation using their infinitesimal generators. Carrying out feedback design based on invariant tracking errors leads to control laws preserving the symmetry group admitted by the considered system. Following this invariant tracking approach proposed in \cite{RR99, MRR04} the present paper contributes to the application of the normalization procedure used to construct invariant tracking errors by giving a geometric interpretation of the so-called moving frame in terms of the group orbits and the meaning of $G$-compatible outputs. It turns out that the normalization procedure can also be applied in order to determine suitable coordinate transformations for a local  reduced order realization of a control system with Lie symmetries.  Finally, the idea of controlled symmetries, i.e.\ ``injecting" symmetries into a given control problem by suitable state feedback, has been illustrated by an example of an invariant control of a bioreactor.


\end{document}